\documentclass[10pt]{amsart}
\usepackage{amsfonts}
\usepackage{amsmath}
\usepackage{amssymb}
\usepackage{color}
\usepackage{epsfig,exscale}
\usepackage[all]{xy}
\usepackage{graphicx}
\usepackage{xspace}

 \newcommand{\sigmat}{\widetilde\sigma}
 \font \eightrm=cmr8

\newtheorem{prop}{Proposition}[section]
\newtheorem{lem}{Lemma} [section]
\newtheorem{cor}{Corollary}[section]

\newtheorem{thm}{Theorem}[section]

\DeclareMathOperator{\End}{End}     

\newcommand{\id}{\mathrm{id}}       

\newcommand{\C}{\mathbb{C}}         
\newcommand{\N}{\mathbb{N}}         
\newcommand{\Q}{\mathbb{Q}}         

\def\<#1,#2>{\langle#1,#2\rangle}   
\newcommand{\sepword}[1]{\quad\mbox{#1}\quad} 

\begin{document}

\title[A noncommutative Bohnenblust--Spitzer identity]
      {A noncommutative Bohnenblust--Spitzer identity for Rota--Baxter algebras
       solves Bogoliubov's recursion}

\author{Kurusch Ebrahimi-Fard}
\address{Max Planck Institute for Mathematics,
         Vivatsgasse 7,
         D-53111 Bonn, Germany.}
\email{kurusch@mpim-bonn.mpg.de}
\urladdr{http://www.th.physik.uni-bonn.de/th/People/fard/}

\author{Dominique Manchon}
\address{Universit\'e Blaise Pascal,
         C.N.R.S.-UMR 6620,
         63177 Aubi\`ere, France}
         \email{manchon@math.univ-bpclermont.fr}
         \urladdr{http://math.univ-bpclermont.fr/~manchon/}

\author{Fr\'ed\'eric Patras}
\address{Laboratoire J.-A. Dieudonn\'e
         UMR 6621, CNRS,
         Parc Valrose,
         06108 Nice Cedex 02, France}
\email{patras@unice.fr} \urladdr{www-math.unice.fr/~patras}

\date{\today}


\begin{abstract}
The Bogoliubov recursion is a particular procedure appearing in
the process of renormalization in perturbative quantum field
theory. It provides convergent expressions for otherwise divergent
integrals. We develop here a theory of functional identities for
noncommutative Rota--Baxter algebras which is shown to encode,
among others, this process in the context of Connes--Kreimer's
Hopf algebra of renormalization. Our results generalize the
seminal Cartier--Rota theory of classical Spitzer-type identities
for commutative Rota--Baxter algebras. In the classical,
commutative, case, these identities can be understood as deriving
from the theory of symmetric functions. Here, we show that an
analogous property holds for noncommutative Rota--Baxter algebras.
That is, we show that functional identities in the noncommutative
setting can be derived from the theory of noncommutative symmetric
functions. Lie idempotents, and particularly the Dynkin idempotent
play a crucial role in the process. Their action on the
pro-unipotent groups such as those of perturbative renormalization
is described in detail along the way.
\end{abstract}

\maketitle



\tableofcontents

\thispagestyle{empty}


\section{Introduction}
\label{sect:Introduction}

Spitzer identities first appeared in fluctuation theory, together
with the notion of Baxter relations ---now called Rota--Baxter
(RB) relations~\cite{spitzer1956,baxter1960,atkinson1963}. It was
soon realized by Rota, Cartier, and others, that the theory could
be founded on purely algebraic grounds and had many other
applications~\cite{rota1969,cartier1972,rotasmith1972}
---appearing retrospectively as one of the many striking successes
of Rota's approach to algebra, combinatorics and their
applications.

The purpose of the present article is to extend the theory to the
noncommutative setting. Indeed, the classical Spitzer identities
involve commutative Rota--Baxter operators and algebras. However,
to consider the noncommutative case is natural. For example, the
integration operator acting on matrix algebras is a RB operator
(in that particular case, the Rota--Baxter relation identifies
with the integration by parts rule), so that the Magnus or
Strichartz identities for the solutions of first order linear
differential equations can be viewed as particular examples of
RB-type identities (see \cite{cefg2006,magnus1954,strich1987} and
our account in section 8 of the present article).

Actually, a striking application of the RB formalism in the
noncommutative setting emerged very recently in the context of the
Connes--Kreimer Hopf algebra approach to renormalization in
perturbative quantum field theory
(pQFT)~\cite{ck2000,egk2004,egk2005,ek2005,manchon2001}, and
motivated the present article. The Bogoliubov recursion is a
purely combinatorial recursive process that allows to give a
meaning to the divergent integrals appearing in
pQFT~\cite{collins1984}. Using the RB point of view, the process
finds a very compact and simple formulation. Abstractly, the
recursion takes place in a particular noncommutative RB algebra
and writes
$$
    X = 1 - R\bigl( X \star a\bigr),
$$
where $X$ is the quantity to be computed recursively, $R$ is the
RB operation, $\star$ is the algebra product, and $a$ is a
divergent series naturally associated to the physical quantities
to be computed (the regularized Feynman rules in, say, the
dimensional regularization scheme \cite{collins1984}). We refer
the reader to the last section of the present article for
definitions and further details.

The above functional identity lies in fact at the heart of
Baxter's original work and, for a commutative RB algebra, its
solution is given by nothing but Spitzer's classical
identity~\cite{spitzer1956}. However, as we mentioned, in
renormalization the very RB structure one has to deal with is
noncommutative, so that these results do not apply.
In~\cite{egk2004} it was shown that one can prove a first
noncommutative Spitzer (also known as Pollaczeck--Spitzer)
identity, and that this identity was related to a so-called
Baker--Campbell--Hausdorff (BCH) recursion, which is another way,
besides Bogoliubov's, to perform recursively the renormalization
process~\cite{egk2005,egm2006}.

Here, we show that one can actually prove more. That is, we derive
noncommutative generalizations of the Bohnenblust--Spitzer
identity (Thm.~\ref{thm:ncBS1}, Thm.~\ref{NCSpitzer},
Thm.~\ref{newNCSpitzer}). They allow us to solve completely the
renormalization problem, in the sense that they lead to a closed
formula for $X$ ---as opposed to the BCH or Bogoliubov
recursions---, and also to the Pollaczeck--Spitzer identity. The
new formula is also completely different from the celebrated
Zimmermann forest formula~\cite{zimmermann1969}, that relies on
particular combinatorial properties of Feynman diagrams.\\

However, although our first motivation was renormalization theory
in the Connes--Kreimer approach, one should be aware that the
existence of noncommutative Bohnenblust--Spitzer identities, as
well as the ideas developed to prove the identities, are of general
interest and should lead to a noncommutative approach in the many
fields where commutative Rota--Baxter algebras have been a useful
tool.

To understand our approach, recall one of the main events in the
study of RB algebras. In their seminal 1972 article
\cite{rotasmith1972}, Rota and Smith showed that Spitzer-type
identities for commutative Rota--Baxter algebras could be
understood as deriving from the theory of symmetric functions.
Here, we actually show that the same is true for noncommutative RB
algebras. That is, we show that functional identities for these
algebras can be derived from the theory of noncommutative
symmetric functions~\cite{gelfand1995} or, equivalently, from the
theory of descent algebras of bialgebras ---a cornerstone of the
modern approach to the theory of free Lie
algebras~\cite{reutenauer1993,patras1994}. In the process, we
establish a connection between noncommutative RB algebras and
quasi-symmetric functions in noncommutative variables.

Eventually, as already alluded at, our findings lead to
noncommutative generalizations of the classical
Bohnenblust--Spitzer identity. Moreover, these new identities are
derived from a functional equation (Thm.~\ref{thm:invDynkinRB})
for the classical Dynkin operator (the iteration of the Lie
bracket in the theory of free Lie algebras). As an application we
present a closed formula for the Bogoliubov recursion in the
context of Connes--Kreimer's Hopf algebra approach to perturbative
renormalization. This last finding is complementary to the main
result of the recent article~\cite{eunomia2006}, where two of the
present authors together with J.M.~Gracia-Bond\'{\i}a proved that
the mathematical properties of locality and the so-called
beta-function in pQFT could be derived from the properties of the
Dynkin operator. Our new findings reinforce the idea that the
Dynkin operator and its algebraic properties have to be considered
as one of the building blocks of the modern mathematical theory of
renormalization.

Some partial results were announced in the Letter~\cite{egp2007}.
We give here their complete proofs and develop the general theory
of noncommutative functional identities for Rota--Baxter algebras
as well as their applications to perturbative renormalization
which were alluded at also in~\cite{egp2007}.\\

Let us briefly outline the organization of this paper. The second
section develops the theory of pro-unipotent groups and
pro-unipotent Lie algebras from the point of view of Lie
idempotents. The properties of the Dynkin idempotent are recalled.
The next section surveys the classical Bogoliubov recursion,
emphasizing the Rota--Baxter approach and the connections with
Atkinson's recursion and Spitzer identities. Definition and basic
properties of Rota--Baxter algebras are recalled in the process,
i.e. the Rota--Baxter double and pre-Lie product. We then extend
in the fourth section Rota--Smith's construction of the free
commutative Rota--Baxter algebra (in an arbitrary number of
generators) to the noncommutative case. This leads naturally to
the link between free Rota--Baxter algebras and noncommutative
symmetric and quasi-symmetric functions that are explored
afterwards. Section 6 concentrates on the Bohnenblust--Spitzer
identity in the noncommutative setting, whereas the following one
features an extension of the Magnus recursion and of Strichartz'
solution thereof to arbitrary noncommutative RB algebras.
Applications to perturbative renormalization are considered in the
last section; the relevance of our general results in this setting
is emphasized, since they provide a non recursive solution to the
computation of counterterms in pQFT together with new theoretical
and computational insights on the subject.


\section{Lie idempotents actions on pro-unipotent groups}
\label{sect:LieIdemp}

The ground field $\mathbb{K}$ over which all algebraic structures
are defined is of characteristic zero.

Lie idempotents are well known to be one of the building blocks of
the modern theory of free Lie algebras. This includes the modern
approaches to the Baker--Campbell--Hausdorff formula, the Dynkin
formula for the Hausdorff series, the Zassenhaus formula, and
continuous versions of the same formulas such as Magnus'
continuous Baker--Campbell--Hausdorff formula~\cite[Chap.
3]{reutenauer1993}.

One of the purposes of the present article is to show that the
same result holds for general Rota--Baxter algebras. That is,
functional identities for Rota--Baxter algebras can be derived
from the theory of Lie idempotents. In a certain sense, the result
is not so surprising: after all, in their seminal work, Rota and
Smith~\cite{rotasmith1972} explained the classical Spitzer
identities for commutative RB algebras by means of the Waring
formula, which holds in the algebra of symmetric functions. Lie
idempotents appear (very generally) as soon as one moves from the
commutative algebra setting to the noncommutative one, where
phenomena such as the Baker--Campbell--Hausdorff formula require,
for their solution and combinatorial expansion, free Lie algebraic
tools. For that purpose, the algebra of symmetric functions, which
encodes many of the main functional identities for commutative
algebras, has to be replaced by the descent algebra, which is the
algebra in which Lie idempotents live naturally. We refer to
Reutenauer's standard reference~\cite{reutenauer1993} for further
details on the subject and a general picture of Lie idempotents,
descent algebras, and their applications to the study of Lie and
noncommutative algebras. The definitions which are necessary for
our purposes are given below.

\medskip

To start with, let us first recall from
\cite{reutenauer1993,gelfand1995,patreu2002} and
\cite{eunomia2006} some definitions and properties relating the
classical Dynkin operator to the fine structure theory of Hopf
algebras. Some details are needed, since the results gathered in
the literature are not necessarily stated in a form convenient for
our purposes.

Recall first the classical definition of the Dynkin operator $D$.
The Dynkin operator is the linear map from $A:=T(X)$, the tensor
algebra over a countable set $X$, into itself defined as the
left-to-right iteration of the associated Lie bracket, so that,
for any sequence $y_1,\dots,y_n$ of elements of $X$:
$$
    D(y_1 \dots y_n):=[\dots[[y_1,y_2],y_3]\dots,y_n]
$$
where $[x,y]:=xy-yx$. We also write $D_n$ for the action of $D$ on
$T_n(X)$, the component of degree $n$ of the tensor algebra (the
linear span of words $y_1 \ldots y_n$, $y_i \in X,\
i=1,\ldots,n$). Notice, for further use, the iterated structure of
the definition of $D$; it will appear below that the Dynkin
operator is a natural object to understand advanced properties of
the Spitzer algebra~\cite{egp2007}. The Dynkin operator can be
shown to be a quasi-idempotent. That is, its action on an
homogeneous element of degree $n$ of the tensor algebra satisfies:
$D^2 = n \cdot D$ and, moreover, the associated projector
$\frac{D}{n}$ is a projection from $T_n(X)$ onto the component of
degree $n$, $Lie_n(X)$, of the free Lie algebra over $X$,
see~\cite{reutenauer1993}.

The tensor algebra is a graded cocommutative connected Hopf
algebra: the coalgebra structure is entirely specified by the
requirement that the elements $x\in X$ are primitive elements in
$T(X)$. It is therefore naturally provided with an antipode $S$
and with a grading operation $Y$ ---the map $Y$ acting as the
multiplication by $n$ on $T_n(X)$. One can then show that the
Dynkin operator can be rewritten in purely Hopf algebraic terms as
$D = S \star Y$, where $\star$ stands for the convolution product
in $\End(T(X))$. Recall, for further use, that, writing $\Delta$
and $\pi$ for the coproduct and the product, respectively, in an
arbitrary Hopf algebra $H$, the convolution product of two
endomorphisms $f$ and $g$ of $H$ is defined by:
$$
    f \star g := \pi \circ (f \otimes g)\circ \Delta.
$$
The definition of the Dynkin map as a convolution product can be
extended to any graded connected cocommutative or commutative Hopf
algebra~\cite{patreu2002}, as a particular case of a more general
phenomenon, namely the possibility to define an action of
$\Sigma_n$, the classical Solomon's algebra of type $A_n$ (resp.
the opposite algebra) on any graded connected commutative (resp.
cocommutative) Hopf algebra~\cite{patras1994}. In particular, if
we call {\textit{descent algebra}} and write
${\mathcal{D}}:=\bigoplus_{n\in\N}{\mathcal{D} }_n$ for the
convolution subalgebra of $\End(T(X))$ generated by the graded
projections, $p_n: T(X) \longrightarrow T_n(X)$, then
$\Sigma_n^{op}$, the opposite algebra to Solomon's algebra of type
$A_n$, identifies naturally with ${\mathcal{D}}_n$, which inherits
an associative algebra structure from the composition product in
$\End(T(X)$ \cite{reutenauer1993,patras1994}.

In the present article we call {\textit{Lie idempotents}} the
projectors from $T_n(X)$ to $Lie_n(X)$ that belong to Solomon's
algebra $\Sigma_n$ (one sometimes calls Lie idempotents the more
general projectors belonging to the symmetric group algebra $\Q
[S_n]$, in which Solomon's algebra is naturally embedded
---however the latter idempotents can not be generalized naturally
to idempotents acting on bialgebras and are therefore not relevant
for our purposes).

A Lie idempotent series is a sequence of Lie idempotents or,
equivalently, a projection map $\gamma$ from $T(X)$ to $Lie(X)$
belonging to the descent algebra the graded components $\gamma_n$
of which are Lie idempotents. Besides the Dynkin idempotent
series, $\frac{D_n}{n}$, the list of Lie idempotents series
include Solomon's Eulerian idempotent series, the Klyachko
idempotent series, the Zassenhaus idempotent series etc. We refer
to \cite{reutenauer1993,gelfand1995} for further details on the
subject.

\begin{prop}\label{propfree}
The descent algebra is a free graded associative algebra freely generated
by the $p_n$. Any Lie idempotent series generates freely the
descent algebra.
\end{prop}
The first part of the Proposition is Corollary~9.14 of
\cite{reutenauer1993}. The second part follows e.g. from
Theorem~5.15 in \cite{gelfand1995}.

As already mentioned, it follows from \cite{patras1994} that all
these idempotent series act naturally on any graded connected
commutative or cocommutative bialgebra. Here, we will restrict our
attention to the cocommutative case, we refer
to~\cite{eunomia2006} for applications of the Dynkin operator
formalism to the commutative but noncocommutative Hopf algebras
appearing in the Connes--Kreimer Hopf algebraic theory of
renormalization in pQFT, see~\cite{ck1998,ck2000,ck2001,ek2005}.

\begin{thm}
\label{thm:invDynkinRB} Let $H=\bigoplus_{n\in\N}H_n$ be an
arbitrary graded connected cocommutative Hopf algebra over a field
of characteristic zero. Any Lie idempotent series induces an
isomorphism between the pro-unipotent group $G(H)$ of group-like
elements of ${\hat{H}}:=\prod_{n\in\N}H_n$ and the pro-nilpotent
Lie algebra $Prim(H)$ of primitive elements in $\hat{H}$.

When the Lie series is the Eulerian idempotent series, the
isomorphism is simply the exponential/logarithm isomorphism
between a pro-unipotent group and its pro-nilpotent Lie algebra.
When the Lie series is the Dynkin series, the inverse morphism is
given by $\Gamma: Prim(H) \rightarrow G(H)$, mapping $h=\sum_{n
\ge 0} h_n$ to:
$$
    \Gamma(h) = \sum\limits_{n \ge 0}
                \sum\limits_{{i_1+ \cdots +i_k = n}\atop i_j>0}
                \frac{h_{i_1} \cdots h_{i_k}}{i_1 (i_1+i_2) \cdots (i_1+\cdots+i_k)}.
$$
\end{thm}

This is a result dual to Theorem~4.1 in~\cite{eunomia2006}, which
established the same formula for characters and infinitesimal
characters of graded connected commutative Hopf algebras when the
Lie series was the Dynkin series. The assertion on pro-unipotency
and pro-nilpotency follows e.g. from the classical equivalence
between group schemes and commutative Hopf algebras. It is simply
a way to recall that the group and Lie algebra we consider come
from graded connected Hopf algebras, and inherit from them the
usual nilpotence and completeness properties of graded connected
algebras (see e.g. \cite{ck1998,ck2000,ck2001,cm2004,cm2006},
where the pro-unipotent group schemes point of view is put to use
systematically instead of the Hopf algebraic one to deal with
similar questions).

We sketch the proof. The descent algebra carries naturally a Hopf
algebra structure \cite{gelfand1995,malreu1995}. Since
$\mathcal{D}$ is freely generated by the $p_n$, the coproduct
$\Delta$ is entirely defined by the requirement that the $p_n$
form a sequence of divided powers (that is, $\Delta
(p_n)=\sum_{i+j=n}p_i\otimes p_j$) or, equivalently, that any Lie
idempotent is a primitive element, see e.g. Corollary~5.17 in
\cite{gelfand1995} or Corollary~3 in \cite{patreu2002}.

It follows from \cite{patreu2002} that there is a compatibility
relation between this coproduct and the descent algebra natural action on
an arbitrary graded connected cocommutative Hopf algebra $H$.
Namely, for any element $f$ in the descent algebra, we have
$$
    \Delta (f) \circ \delta = \delta \circ f
$$
where $\delta$ stands for the coproduct in $H$, and where the
action of $f$ on $H$ is induced by the convolution algebra
morphism that maps $p_n$, viewed as an element of $\mathcal D$, to
the graded projection (also written abusively $p_n$) from $H$ to
$H_n$. In particular, for any Lie idempotent $l_n$ acting on $H_n$ and
any $h \in H_n$, we have:
 \allowdisplaybreaks{
\begin{eqnarray*}
    \delta (l_n(h))&=&\Delta (l_n)(\delta (h))\\
                   &=&(l_n \otimes \epsilon + \epsilon \otimes l_n)
                      (h \otimes 1 + 1 \otimes h + h'\otimes h'')\\
                   &=& l_n(h) \otimes 1+1\otimes l_n(h)
\end{eqnarray*}}
where $\epsilon$ stands for the counit of $H$ (the natural
projection from $H$ to $H_0$ with kernel
$H^+:=\bigoplus_{n>0}H_n$) and $h' \otimes h''$ belongs to $H^+
\otimes H^+$; the identity follows from $l_n$ being primitive in
the descent algebra. So that, in particular, $l_n(h) \in Prim(H)$.
This implies that the action of any Lie series on $H$ and, in
particular, on $G(H)$, the set of group-like elements in $H$,
induces a map to $Lie(H)$.

The particular example of the Eulerian idempotent is interesting.
The Eulerian idempotent $e$ is the logarithm of the identity in
the endomorphism algebra of $T(X)$ and belongs to the descent
algebra. It acts on $H$ as the logarithm of the identity of $H$ in
the convolution algebra $End(H)$, see
\cite{solomon1968,reutenauer1993,patras1993,patras1994}. For any
$h \in G(H)$, we get:
\begin{eqnarray*}
    e(h)=\log(Id_H)(h)&=&\sum\limits_{n\in\N}\frac{(-1)^{n-1}}{n}(Id_H-\epsilon)^{\star n}(h)
                      =\sum\limits_{n\in\N}\frac{(-1)^{n-1}}{n}(h-\epsilon (h))^n=\log(h),
\end{eqnarray*}
since, $h$ being group-like: $(f \star g)(h)=f(h)g(h)$ for any
$f,g \in End(H)$, so that $(Id_H - \epsilon)^{\star
n}(h)=(h-\epsilon (h))^n$. This proves the assertion on the
Eulerian idempotent in the Theorem~\ref{thm:invDynkinRB}.

Now, according to Proposition~\ref{propfree}, the graded components of
any Lie series $l$ generate freely the descent algebra as an
associative algebra. It follows in particular that the identity of
$T(X)$ can be written as a noncommutative polynomial in the $l_n$.
That is, for any Lie series $l$ there exist a unique family of
coefficients $\alpha_{n_1,\ldots,n_k}$ such that:
$$
    Id_{T(X)}=\sum\limits_{n=0}^\infty
              \sum\limits_{n_1+\cdots+n_k=n} \alpha_{n_1,\dots,n_k}l_{n_1}\star \cdots \star l_{n_k}.
$$
In particular, for any $g \in G(H)$, we get, since the element
$Id_{T(X)}$ of the descent algebra acts as the identity on $H$
\cite{patras1994}:
$$
    g=\sum\limits_{n=0}^\infty
      \sum\limits_{n_1+\dots+n_k=n}\alpha_{n_1,\dots,n_k}l_{n_1}(g)\star \cdots \star l_{n_k}(g)
$$
so that the map from $Lie(H)$ to $G(H)$:
$$
    Lie(H)=\bigoplus\limits_{n=0}^\infty Lie_n(H) \ni \sum\limits_{n=0}^\infty
    \lambda_n \longmapsto \sum\limits_{n=0}^\infty
                          \sum\limits_{n_1+\dots+n_k=n}
                          \alpha_{n_1,\dots,n_k}\lambda_{n_1}\star \cdots \star \lambda_{n_k}
$$
is a right inverse (and in fact also a left inverse) to $l$. The
particular formula for the Dynkin idempotent follows from
\cite{gelfand1995} or from Lemma~2.1 in \cite{eunomia2006}.

\medskip

Two particular applications of the theorem are well-known.
Consider first the case where $H$ is the Hopf algebra of
noncommutative symmetric functions. Then, $H$ is generated as a
free associative algebra by the {\textit{complete homogeneous
noncommutative symmetric functions}} $S_k$, $k \in \N$, which form
a sequence of divided powers, that is, their sum is a group-like
element in $\hat{H}$. The graded components of the corresponding
primitive elements under the action of the Dynkin operator are
known as the {\textit{power sums noncommutative symmetric
functions of the first kind}} \cite{gelfand1995}. Second, consider
the classical descent algebra viewed as a Hopf algebra. Then, the
Dynkin operator (viewed as the convolution product $S\star Y$
acting on the Hopf algebra $\mathcal D$) sends the identity of
$T(X)$, which is a group-like element of the descent algebra to
the classical Dynkin operator. This property was put to use in
Reutenauer's monograph to rederive all the classical functional
Lie-type identities in the tensor algebra ---for example the
various identities related to the Baker--Campbell--Hausdorff
formula.

As it will appear, a surprising conclusion of the present article
is that {\textit{the same machinery}} can be used to derive the
already known formulas for commutative Rota--Baxter algebras but,
moreover, can be used to prove new formulas in the noncommutative
setting.


\section{Rota--Baxter algebras and Bogoliubov's recursion}
\label{sect:RB}

In this section we first briefly recall the definition of
Rota--Baxter (RB) algebra and its most important properties. For
more details we refer the reader to the classical
papers~\cite{atkinson1963,baxter1960,cartier1972,rota1969,rotasmith1972},
as well as for instance to the references~\cite{ek2005,egm2006}.

Let $A$ be an associative not necessarily unital nor commutative
algebra with $R \in \End(A)$. The product of $a$ and $b$ in $A$ is
written $a \cdot b$ or simply $ab$ when no confusion can arise. We
call a tuple $(A,R)$ a {\textit{Rota--Baxter algebra}} of weight
$\theta \in \mathbb{K}$ if $R$ satisfies the {\textit{Rota--Baxter
relation}}
\begin{equation}
    R(x)R(y) = R\bigl( R(x) y + x R(y) + \theta xy \bigr).
\label{def:RBR}
\end{equation}
Changing $R$ to $R' := \mu R$, $\mu \in \mathbb{K}$, gives rise to
a RB algebra of weight $\theta':=\mu \theta$, so that a change in
the $\theta$ parameter can always be achieved, at least as long as
weight non-zero RB algebras are considered.

The definition generalizes to other types of algebras than
associative algebras: for example one may want to consider RB Lie
algebra structures. Further below we will encounter examples of
such structures.

In the following we denote the particular argument of the map $R$
on the right hand side of (\ref{def:RBR}) by
$$
    x \ast_\theta y := R(x) y + x R(y) + \theta xy,
$$
and will come back to it further below.

Let us recall some classical examples of RB algebras. First,
consider the integration by parts rule for the Riemann integral
map. Let $A := C(\mathbb{R})$ be the ring of real continuous
functions with pointwise product. The indefinite Riemann integral
can be seen as a linear map on~$A$
\begin{equation}
    I: A \to A, \qquad\  I(f)(x) := \int_0^x f(t)\,dt.
\label{eq:Riemann}
\end{equation}
Then, integration by parts for the Riemann integral can be written
compactly as
\begin{equation}
    I(f)(x)I(g)(x) = I\bigl( I(f) g \bigr)(x) + I\bigl( fI(g)\bigr)(x),
\label{eq:integ-by-parts}
\end{equation}
dually to the classical Leibniz rule for derivations. Hence, we
found our first example of a weight zero Rota--Baxter map.
Correspondingly, on a suitable class of functions, we define the
following Riemann summation operators
\begin{eqnarray}
    R_\theta(f)(x) := \sum_{n = 1}^{[x/\theta]} \theta f(n\theta)
    \qquad\ {\rm{and}} \qquad\
    R'_\theta(f)(x) := \sum_{n = 1}^{[x/\theta]-1} \theta f(n\theta).
\label{eq:le-clou}
\end{eqnarray}
We observe readily that
\begin{align}
    &\biggl( \sum_{n = 1}^{[x/\theta]} \theta f(n\theta)\biggr)
     \biggl( \sum_{m = 1}^{[x/\theta]} \theta g(m\theta)\biggr)
   = \biggl( \sum_{n > m = 1}^{[x/\theta]}
        + \sum_{m > n = 1}^{[x/\theta]}
        + \sum_{m = n = 1}^{[x/\theta]}  \biggr)\theta^2 f(n\theta) g(m\theta) \nonumber \\
    &= \sum_{m = 1}^{[x/\theta]} \theta^2 \biggl(\sum_{k = 1}^{m} f\bigl(k\theta\bigr)\biggr)
                                                        g(m\theta)
     + \sum_{n = 1}^{[x/\theta]} \theta^2 \biggl(\sum_{k = 1}^{n} g\bigl(k\theta\bigr)\biggr)
                                                                             f(n\theta) \nonumber
     - \sum_{n = 1}^{[x/\theta]} \theta^2 f(n\theta)g(n\theta)\\
    &= R_\theta\bigl(R_\theta(f)g\bigr)(x) + R_\theta\bigl(fR_\theta(g)\bigr)(x) + \theta R_\theta(fg)(x).
\label{Riemsum1}
\end{align}
Similarly for the map $R'_\theta$. Hence, the Riemann summation
maps $R_\theta$ and $R'_\theta$ satisfy the weight $-\theta$ and
the weight $\theta$ Rota--Baxter relation, respectively.

Another classical example, and the reason why RB algebras first
appeared in fluctuation theory, comes from the operation that
associates to the characteristic function of a real valued random
variable $X$ the characteristic function of the random variable
$max(0,X)$. It is worth pointing out that all these classical
examples involve {\textit{commutative}} RB algebras.

One readily verifies that $\tilde{R}:=- \theta \id_A - R$ is a
Rota--Baxter operator. Note that
$$
    R(a) \tilde{R}(b) = \tilde{R}(R(a)b) + R(a\tilde{R}(b)),
$$
and similarly exchanging $R$ and~$\tilde{R}$. In the following we
denote the image of $R$ and $\tilde{R}$ by $A_{-}$ and $A_+$,
respectively.

\begin{prop}
Let $(A,R)$ be a Rota--Baxter algebra. $A_{\pm} \subseteq A$ are
subalgebras in $A$.
\end{prop}

We omit the proof since it follows directly from the Rota--Baxter
relation. A {\textit{Rota--Baxter ideal}} of a Rota--Baxter
algebra $(A,R)$ is an ideal $I \subset A$ such that $R(I)
\subseteq I$.

The Rota--Baxter relation extends to the Lie algebra $L_A$
corresponding to $A$
$$
    [R(x), R(y)] = R\bigl([R(x),y] + [x,R(y)]\bigr) + \theta R\bigl([x, y]\bigr)
$$
making $(L_A,R)$ into a Rota--Baxter Lie algebra of weight
$\theta$. Let us come back to the product we defined after
equation~(\ref{def:RBR}).

\begin{prop}
The vector space underlying $A$ equipped with the product
\begin{eqnarray}
        x \ast_\theta y := R(x)y + xR(y) + \theta xy
\label{def:RBdouble}
\end{eqnarray}
is again a Rota--Baxter algebra of weight $\theta$ with
Rota--Baxter map $R$. We denote it by $(A_\theta,R)$ and call it
{\textit{double Rota--Baxter algebra}}.
\end{prop}

\begin{proof}
Let $x,y,z \in A$. We first show
associativity
 \allowdisplaybreaks{
\begin{eqnarray}
    x \ast_\theta ( y \ast_\theta z)
    &=& x R\bigl( yR(z) + R(y)z + \theta yz \bigr) + R(x) \bigl( yR(z) + R(y)z + \theta yz \bigr) \nonumber\\
    & & \hspace{1cm} + \theta x \bigl( yR(z) + R(y)z + \theta yz \bigr)    \nonumber\\
    &=& xR(y)R(z) + R(x)yR(z) + \theta xyR(z) + R(x)R(y)z                  \nonumber\\
    & & \hspace{1cm} + \theta xR(y)z + \theta R(x)yz + \theta^2 xyz        \nonumber\\
    &=& (x \ast_\theta y) \ast_\theta z.                                   \nonumber
\end{eqnarray}}
Now we show that the original Rota--Baxter map $R \in \End(A)$
also fulfills the Rota--Baxter relation with respect to the
$\ast_\theta$-product.
 \allowdisplaybreaks{
\begin{eqnarray*}
    R(x) \ast_\theta R(y) - \theta R(x \ast_\theta  y)
    &=& R^2(x)R(y) + R(x)R^2(y) + \theta R(x)R(y) - \theta R(x)R(y)\\
    &=& R\bigl( x \ast_\theta R(y) + R(x) \ast_\theta y \bigr).
\end{eqnarray*}}
\end{proof}
We used the following homomorphism property of the Rota--Baxter
map between the algebras $A_{\theta}$ and $A$.

\begin{lem} \label{lem:RBhomomorph}
Let $(A,R)$ be a Rota--Baxter algebra of weight $\theta$. The
Rota--Baxter map $R$ becomes a (not necessarily unital even if $A$
is unital) algebra homomorphism from the algebra $A_{\theta}$ to
$A$
 \allowdisplaybreaks{
\begin{eqnarray}
    R\bigl(a \ast_\theta b\bigr) &=& R(a)R(b).
    \label{eq:RBhom1}
\end{eqnarray}}
For $\tilde{R}:= - \theta \id_A - R$ we find
\begin{equation}
  \tilde{R}\bigl( a \ast_\theta b \bigr) = -\tilde{R}(a)\tilde{R}(b).
  \label{eq:RBhom2}
\end{equation}
\end{lem}
We remark here that if $R$ is supposed to be idempotent, $(A,R)$
must be a Rota--Baxter algebra of unital weight $\theta=-1$. Now
we introduce the notion of an associator for $x,y,z \in B$, where
$B$ is an arbitrary algebra:
$$
    a_{\cdot}(x,y,z):=(x \cdot y) \cdot z - x \cdot (y \cdot z).
$$
Recall that a left {\textit{pre-Lie algebra}} $P$ is a vector
space, together with a bilinear pre-Lie product $\triangleright: P
\otimes P \to P$, satisfying the left pre-Lie relation
$$
        a_{\triangleright}(x,y,z)=a_{\triangleright}(y,x,z).
$$
With an obvious analog notion of a right pre-Lie product
$\triangleleft: P \otimes P \to P$ and right pre-Lie relation
$$
        a_{\triangleleft}(x,y,z)=a_{\triangleleft}(x,z,y).
$$
See~\cite{chapoliv2001} for more details. Let $P$ be a left (or
right) pre-Lie algebra. The commutator $[a,b]_\triangleright:= a
\triangleright b - b \triangleright a$ for $a,b \in P$ satisfies
the Jacobi identity. Hence, the vector space $P$ together with
this commutator is a Lie algebra, denoted by $L_P$. Of course,
every associative algebra is also pre-Lie.

\begin{lem} \label{lem:pre-LieRB}
Let $(A,R)$ be an associative Rota--Baxter algebra of weight
$\theta$. The binary compositions
 \allowdisplaybreaks{
\begin{eqnarray}
  a \triangleright_\theta b &:=& R(a)b - bR(a) - \theta ba
                             = [R(a),b] - \theta ba
                            = R(a)b - \bigl(- b\tilde{R}(a)\bigr),
  \label{lem:leftRBpreLie}\\
  a \triangleleft_\theta b &:=& aR(b) - R(b)a - \theta ba
                      = [a,R(b)] - \theta ba
                      = aR(b) - \bigl(- \tilde{R}(b)a\bigr),
  \label{lem:rightRBpreLie}
\end{eqnarray}}
define a left respectively right pre-Lie structure on $A$.
\end{lem}

\begin{proof} The Lemma follows by direct inspection. It may also be
deduced from deeper properties of Rota--Baxter algebras related to
dendriform di- and
trialgebras~\cite{aguiar2000,e2002,loday2001,lodayronco2002}. That
is, the identification $a \prec b:= aR(b)$, $a \bullet b = \theta
ab$ and $a \succ b:=R(a)b$ in a Rota--Baxter algebra $(B,R)$
defines a dendriform trialgebra, hence also a dialgebra structure.
We refer the reader to~\cite{KDF2007} for more.
\end{proof}

Recall that anti-symmetrization of a pre-Lie product gives a Lie
bracket. In the case of the Rota--Baxter pre-Lie compositions
(\ref{lem:leftRBpreLie},\ref{lem:rightRBpreLie}), we find
 \allowdisplaybreaks{
\begin{eqnarray*}
    [a,b]_{\triangleright_\theta} &:=& a \triangleright_\theta b  - b \triangleright_\theta a
                                   =  a \triangleleft_\theta b  - b \triangleleft_\theta a
                                   =[a,b]_{\triangleleft_\theta}\\
          &=& [R(a),b] + [a,R(b)] + \theta [a,b] = [a,b]_{\ast_\theta}.
\label{def:RBLiedouble}
\end{eqnarray*}}
Hence, the double Rota--Baxter product and the left as well as
right Rota--Baxter pre-Lie products define the same Lie bracket on
$(A,R)$.

\begin{lem}
Let $(A,R)$ be an associative Rota--Baxter algebra of weight
$\theta$. The left pre-Lie algebra $(A,\triangleright_\theta)$ is
a Rota--Baxter left pre-Lie algebra of weight $\theta$, with
Rota--Baxter map $R$. Similarly for $(A,\triangleleft_\theta)$
being a Rota--Baxter right pre-Lie algebra of weight $\theta$.
\end{lem}

\begin{proof} We prove only the statement for the left RB pre-Lie algebra. Let $x,y \in A$.
 \allowdisplaybreaks{
  \begin{eqnarray*}
        R(x) \triangleright_\theta R(y) &=& R\bigl(R(x)\bigr)R(y) - R(y)R\bigl(R(x)\bigr) - \theta R(y)R(x) \nonumber\\
                           &=& R\bigl( R(R(x))y + R(x)R(y) + \theta R(x)y \bigr)                           \\
                           & & - R\bigl( yR\bigl(R(x)\bigr) + R(y)R(x) + \theta yR(x)\bigr)
                               - \theta R\bigl( R(y)x + yR(x) + \theta yx\bigr)                            \\
                           &=& R\bigl( R(x)\triangleright_\theta y + x \triangleright_\theta R(y)
                               + \theta x \triangleright_\theta y\bigr).
  \end{eqnarray*}}
\end{proof}

Let us now turn to the Bogoliubov recursion. Briefly, this
recursion provides an elaborate procedure to give a sense (or to
renormalize, that is, to associate a finite quantity, called the
renormalized amplitude) to divergent integrals appearing in
perturbative high-energy physics calculations. The renormalization
process and in particular the Bogoliubov recursion can be
reformulated in purely algebraic terms inside the Connes--Kreimer
paradigmatic Hopf algebraic approach to perturbative
renormalization. We follow this point of view and refer the reader
to Collins' monograph \cite{collins1984} and the by now standard
references~\cite{ck1998,ck2000,ck2001,fg2005,manchon2001} for
further information on the subject. Further details on the
physical meaning of the recursion will be given in the last
section of the article, we concentrate for the time being on its
mathematical significance.

Let us outlay the general setting, following~\cite{egp2007}. Let
$H$ be a graded connected commutative Hopf algebra (in the
physical setting $H$ would be a Hopf algebra of Feynman diagrams),
and let $A$ be a commutative unital algebra. Assume further that
$A$ splits into a direct sum of subalgebras, $A = A_+ \bigoplus
A_-$, with $1 \in A_+$. The projectors to $A_\pm$ are written
$R_\pm$ respectively. The pair $(A,R_-)$ is then a weight
$\theta=-1$ commutative Rota--Baxter algebra, whereas the algebra
$Lin(H,A)$ with the idempotent operator defined by ${\mathcal
R}_-(f) := R_- \circ f$ for $f \in Lin(H,A)$ is a (in general
noncommutative) unital Rota--Baxter algebra. Here, the algebra
structure on $Lin(H,A)$ is induced by the convolution product,
that is, for any $f,g \in Lin(H,A)$ and any $h \in H$:
$$
    f\star g(h):=f(h^{(1)})g(h^{(2)}),
$$
with Sweedler's notation for the action of the coproduct $\delta$
of $H$ on $h$: $\delta (h)=h^{(1)}\otimes h^{(2)}$.

The essence of renormalization is contained in the existence of a
decomposition of the group $G(A)$ of algebra maps from $H$ to $A$
into a (set theoretic) product of the groups $G_-(A)$ and $G_+(A)$
of algebra maps from $H^+$ to $A_-$, respectively from $H$ to
$A_+$. We view $G_-(A)$ as a subgroup of $G(A)$ by extending maps
$\gamma$ from $H^+$ to $A_-$ to maps from $H$ to $A$ by requiring
that $\gamma(1)=1$. In concrete terms, any element $\gamma$ of
$G(A)$ can be rewritten uniquely as a product $\gamma_-^{-1} \star
\gamma_+$, where $\gamma_- \in G_-(A)$ and $\gamma_+ \in G_+(A)$.

The Bogoliubov recursion is a process allowing the inductive
construction of the elements $\gamma_-$ and $\gamma_+$ of the
aforementioned decomposition. Writing $e_A$ for the unit of $G(A)$
(the counit map of $H$ composed with the unit map of
$A$, $e_A=\eta_A \circ \epsilon$), $\gamma_-$ and $\gamma_+$ solve
the equations:
\begin{equation}
    \gamma_\pm = e_A \pm {\mathcal R}_\pm \bigl(\gamma_- \star (\gamma -e_A)\bigr).
\label{bogorecursions}
\end{equation}
The recursion process is by induction on the degree $n$ components
of $\gamma_-$ and $\gamma_+$ viewed as elements of the (suitably
completed) graded algebra $Lin(H,A)$. The fact that the recursion
defines elements of $G_-(A)$ and $G_+(A)$ is not obvious from the
definition: since the recursion takes place in $Lin(H,A)$, one
would expect $\gamma_- - e_A$ and $\gamma_+$ to belong to
$Lin(H,A_-)$ and $Lin(H,A_+)$, respectively. The fact that $\gamma_-$ and
$\gamma_+$ do belong to $G_-(A)$ and $G_+(A)$, respectively,
follows from the Rota--Baxter algebra structure of $A$, as has
been shown by Kreimer and Connes-Kreimer, see e.g.
\cite{kreimer1999,ck2000} and the references therein. The map
\begin{equation}
    \bar{\gamma}:=\gamma_- \star (\gamma - e_A)
\label{BogoRbarMap}
\end{equation}
is called {\textit{Bogoliubov's preparation}} or
$\bar{R}$-{\textit{operation}}. Hence, on $H^+$ we see that
$\gamma_{\pm} = \pm {\mathcal R}_\pm (\bar{\gamma})$.

Setting $a:=-(\gamma - e_A)$, the recursion can be rewritten:
\begin{equation}
    \gamma_\pm = e_A \mp {\mathcal R}_\pm(\gamma_- \star a),
\label{BogoRec}
\end{equation}
and can be viewed as an instance of results due to F.V.~Atkinson,
who, following Baxter's work~\cite{baxter1960}, has made an
important observation~\cite{atkinson1963} when he found a
multiplicative decomposition for associative unital Rota--Baxter
algebras. We will state his result for the ring of power
series~$B[[t]]$, $(B,R)$ an arbitrary Rota--Baxter algebra.
Inductively define in a general RB algebra $(B,R)$,
\begin{equation}
    (Ra)^{[n+1]} := R\bigl((Ra)^{[n]}a\bigr) \ \sepword{and}\ (Ra)^{\{n+1\}} := R\bigl(a(Ra)^{\{n\}}\bigr).
\label{connvention}
\end{equation}
with the convention that $(Ra)^{[1]}:=R(a)=:(Ra)^{\{1\}}$ and
$(Ra)^{[0]}:=1_B=:(Ra)^{\{0\}}$.

\begin{thm}
\label{Atkinson1} Let $(B,R)$ be an associative unital
Rota--Baxter algebra. Fix $a \in B$ and let $F$ and~$G$ be defined
by $F:=\sum_{n\in\N} t^n (Ra)^{[n]}$ and $G:=\sum_{n\in\N}
t^n(\tilde{R}a)^{\{n\}}$. Then, they solve the equations
\begin{equation}
    F = 1_{B} + tR(F \ a)
    \sepword{and}
    G = 1_{B} + t\tilde{R}(a \ G),
\label{atkinsonEqs}
\end{equation}
in~$B[[t]]$ and we have the following factorization
\begin{equation}
    F \bigl(1_{B} + at\theta \bigr) G = 1_{B},
     \ \sepword{ so that } \
    1_B + a t \theta = F^{-1}G^{-1}.
\label{eq:Atkinsonfact1}
\end{equation}
For an idempotent Rota--Baxter map this factorization is unique.
\end{thm}

\begin{proof}
The proof follows simply from calculating the product $FG$.
Uniqueness for idempotent Rota--Baxter maps is easy to show, see
for instance~\cite{egm2006}.
\end{proof}

One may well ask what equations are solved by the inverses
$F^{-1}$ and $G^{-1}$. We answer this question, the solution of
which will be important in forthcoming developments, in the
following corollary.

\begin{cor}
\label{inverseAtkinson} Let $(B,R)$ be an associative unital
Rota--Baxter algebra. Fix $a \in B$ and assume~$F$ and~$G$ to
solve the equations in the foregoing theorem. The
inverses~$F^{-1}$ and~$G^{-1}$ solve the equations
\begin{equation}
    F^{-1} = 1_{B} - tR(a \ G)
    \sepword{and}
    G^{-1} = 1_{B} - t\tilde{R}(F \ a),
\label{atkinsonInvEqs}
\end{equation}
in~$B[[t]]$.
\end{cor}

\begin{proof}
Let us check this for $F$ and $F^{-1}$. Recall that $G = 1_{B} +
t\tilde{R}(a \ G)$.
 \allowdisplaybreaks{
  \begin{eqnarray*}
    FF^{-1} &=& 1_{B} - tR(a \ G) + tR(F \ a) - t^2R(F \ a)R(a \ G)\\
            &=& 1_{B} - tR(a \ G) + tR(F \ a)
                      - t^2R\bigl(R(F \ a)a \ G\bigr) -t^2R\bigl(F \ aR(a \ G)\bigr)
                      - t^2\theta R\bigl(F \ a^2 \ G\bigr)\\
            &=& 1_{B} - tR\bigl( (1_B + tR(F \ a))a \ G\bigr)
                      + tR(F \ a)
                      + t^2R\bigl(F \ a\tilde{R}(a \ G)\bigr)\\
            &=& 1_{B} - tR\bigl( (1_B + tR(F \ a))a \ G\bigr)
                  + tR\bigl( F \ a (1_B + t\tilde{R}(a \ G))\bigr)\\
            &=& 1_{B} + tR\bigl( F \ a \ G\bigr) - tR\bigl( F\ a \ G\bigr) =  1_{B}
  \end{eqnarray*}}
\end{proof}

Going back to Bogoliubov's recursions~(\ref{bogorecursions}) we
see that $\gamma_-$ corresponds to the first equation in
(\ref{atkinsonEqs}), whereas $\gamma_+$ corresponds to the inverse
of the second equation in (\ref{atkinsonEqs}), see
(\ref{atkinsonInvEqs}).

The solution to Atkinson's recursion can be simply expressed as
follows: the coefficient of $t^n$ in the expansion of $F$ is
$(Ra)^{[n]}$. When the Rota--Baxter algebra $(A,R)$ is
commutative, the classical Spitzer formulas allow to reexpress and
expand these terms, giving rise to non-recursive expansions. The
first Spitzer identity, or Pollaczeck--Spitzer identity reads:
$$
    \sum\limits_{m\in\N} t^m(Ra)^{[m]} = \exp\bigl(\theta^{-1} R\log (1+at\theta)\bigr),
$$
where $a$ is an arbitrary element in a weight $\theta$
Rota--Baxter algebra $A$~\cite{baxter1960,spitzer1956}. In the
framework of the Rota--Smith presentation~\cite{rotasmith1972} of
the free commutative RB algebra on one generator (the ``standard''
RB algebra), this becomes the Waring formula relating elementary
and power sum symmetric functions~\cite{sagan2001}. In fact,
comparing the coefficient of $t^n$ on both sides, Spitzer's
identity says that:
$$
    n!(Ra)^{[n]} = \sum\limits_\sigma (-\theta)^{n-k(\sigma )} R(a^{|\tau_1|})\cdots R(a^{|\tau_{k(\sigma )}|}),
$$
where the sum is over all permutations on $[n]$ and $\sigma =
\tau_1 \dots \tau_{k(\sigma )}$ is the decomposition of $\sigma$
into disjoint cycles~\cite{rotasmith1972}. We denote by $|\tau_i|$
the number of elements in $\tau_i$. The second Spitzer formula, or
{\textit{Bohnenblust--Spitzer formula}}, follows by polarization
\cite{rotasmith1972}:
$$
    \sum\limits_\sigma R\bigl(R(\dots (R(a_{\sigma_1})a_{\sigma_2}\dots )a_{\sigma_n})\bigr)
    =\sum\limits_{\pi \in{ \mathcal P}_n}(-\theta)^{n-|\pi |}
     \prod\limits_{\pi_i\in \pi}(m_i -1)!R(\prod\limits_{j\in\pi_i}a_j)
$$
for an arbitrary sequence of elements $a_1,\dots,a_n$ in $A$.
Here, $\pi$ runs over unordered set partitions ${ \mathcal P}_n$
of $[n]$; by $|\pi|$ we denote the number of blocks in $\pi$; and
$m_i:=|\pi_i|$ is the size of the particular block $\pi_i$.

In the sequel of the article, we will show how these identities
can be generalized to arbitrary (i.e. noncommutative) Rota--Baxter
algebras, giving rise to closed formulas for the terms in the
Bogoliubov, i.e. Atkinson recursion.


\section{Free Rota--Baxter algebras and NCQSym}
\label{sect:freeNCRB}

In the present section, we introduce a model for noncommutative
(NC) weight one free RB algebras that extends to the
noncommutative setting the notion of {\it{standard Baxter
algebra}}~\cite{rotasmith1972}.

Let $X=(x_1,\dots,x_n,\dots)$ be an ordered set of variables, or
{\textit{alphabet}} and $T(X)$ be once again the tensor algebra or free
associative algebra over $X$. Recall that the elements of $T(X)$ are linear
combinations of noncommutative products $x_{i_1} \dots x_{i_k}$ of
elements of $X$, or words over $X$. We shall also consider finite
ordered families of alphabets $X^1,\dots ,X^n$ and write $x_n^i$
for the elements in $X^i=(x_1^i,\dots,x_n^i,\dots )$. The tensor
algebra over $X^1 \coprod \dots \coprod X^n$ is written
$T(X^1,\dots,X^n)$.

We write $A$ (resp. $A^{(n)}$) for the algebra of countable
sequences $Y=(y_1,\dots ,y_n,\dots )$ of elements of $T(X)$ (resp.
$T(X^1,\dots,X^n)$) equipped with pointwise addition and products:
$(y_1,\dots ,y_n,\dots )+ (z_1,\dots ,z_n,\dots ):=(y_1+z_1,\dots
,y_n+z_n,\dots )$ and $(y_1,\dots ,y_n,\dots ) \cdot (z_1,\dots
,z_n,\dots ):=(y_1 \cdot z_1, \dots , y_n \cdot z_n,\dots )$. We
also write $Y_i$ for the $i$-th component of the sequence $Y$. By
a slight abuse of notation, we view $X$ (resp. $X^i,\ i\leq n$) as
a sequence, and therefore also as an element of $A$ (resp. of
$A^{(n)}$).

\begin{lem}
The operator $R \in \End(A)$ (resp. $R \in \End(A^{(n)})$)
$$
    R(y_1,\dots ,y_n,\dots ):=(0,y_1,y_1+y_2,\dots ,y_1 + \dots  + y_n,\dots )
$$
defines a weight one RB algebra structure on $A$ (resp. $A^{(n)}$).
\end{lem}

\begin{proof}
Let us check the formula ---in the sequel we will omit some
analogous straightforward verifications.
 \allowdisplaybreaks{
\begin{eqnarray*}
  \lefteqn{R\bigl( (y_1, \dots  ,y_n, \dots  ) \cdot R(z_1, \dots  ,z_n, \dots  )\bigr)
     +R\bigl(R(y_1, \dots  , y_n, \dots  ) \cdot (z_1, \dots  ,z_n, \dots  )\bigr)}\\
    & = &\bigl(0,0,y_2z_1, \dots  , \sum \limits_{i=1}^{n-1} y_{i}(z_1 + \dots  + z_{i-1}), \dots  \bigr)
       + \bigl(0,0,y_1z_2, \dots ,\sum\limits_{i=1}^{n-1}(y_1 + \dots  + y_{i-1})z_{i}, \dots  \bigr)\\
    & = &\bigl(0,y_1z_1-y_1z_1,(y_1+y_2)(z_1+z_2)-(y_1z_1+y_2z_2), \dots  \\
    & &\hspace{2cm} \dots , (y_1 + \dots  +y_{n-1})(z_1+ \dots  +z_{n-1})
                              - (y_1z_1+ \dots  +y_{n-1}z_{n-1}),\dots \bigr)\\
    & = & R(y_1, \dots  ,y_n, \dots  ) \cdot R(z_1, \dots  ,z_n, \dots )
       -R\bigl((y_1, \dots  ,y_n, \dots  )\cdot (z_1, \dots  ,z_n, \dots  )\bigr)
\end{eqnarray*}}
\end{proof}

Recall the notation introduced in (\ref{connvention}), where we
defined inductively $(R a)^{[n]}$ and $(R a)^{\{n\}}$. Let us
recall also Hivert's notion of quasi-symmetric functions over a
set of noncommutative variables from
\cite{bergeron2005,novelli2006}. Let $f$ be a surjective map from
$[n]$ to $[k]$, and let $X$ be a countable set of variables, as
above. Then, the quasi-symmetric function over $X$ associated to
$f$, written $M_f$ is, by definition,
$$
    M_f := \sum\limits_\phi x_{\phi^{-1} \circ f(1)} \dots x_{\phi^{-1}\circ f(n)},
$$
where $\phi$ runs over the set of increasing bijections between
subsets of $\N$ of cardinality $k$ and $[k]$. It is often
convenient to represent $f$ as the sequence of its values,
$f=(f(1),\dots ,f(n))$ or $f=f(1),\dots ,f(n)$ in the notation
$M_f$. The definition is best understood by means of an example:
$$
    M_{1,3,3,2}=x_1x_3x_3x_2 + x_1x_4x_4x_2 + x_1x_4x_4x_3 + x_2x_4x_4x_3 + \dots
$$
The linear span $NCQSym$ of the $M_f$'s is a subalgebra of the
algebra of noncommutative polynomials over $X$ (up to classical
completion arguments that we omit, and that allow to deal with
infinite series such as the $M_f$ as if they were usual
noncommutative polynomials). It is related to various fundamental
objects such as the Coxeter complex of type $A_n$ or the
corresponding Solomon--Tits and twisted descent algebras. We refer
to \cite{ps2005,bergeron2005,novelli2006} for further details on
the subject.

We also introduce, for further use, the notation, $M_f^{n}$ for
the image of $M_f$ under the map sending $x_i$ to $0$ for $i>n$
and $x_i$ to itself else. For the above example for instance we
find
$$
    M^3_{1,3,3,2} = x_1x_3x_3x_2,
$$
$$
M^4_{1,3,3,2} = x_1x_3x_3x_2+x_1x_4x_4x_2+x_1x_4x_4x_3+x_2x_4x_4x_3
$$
At last, we write $[n]$ for the identity map on $[n]$ and
$\omega_n$ for the endofunction of $[n]$ reversing the ordering,
so that $\omega_n(i):= n-i+1$ and
$M_{\omega_n}=M_{n,n-1,\dots,1}$.

\begin{prop}
In the RB algebra $A$, we have:
$$
    (RX)^{[n]}=(0,M_{[n]}^{1},M_{[n]}^{2},\dots ,M_{[n]}^{k-1},\dots ),\ n>1,
$$
where $M_{[n]}^{k-1}$ is at the $k$th position in the sequence, and
$$
    (RX)^{\{n\}}=(0,M_{\omega_n}^{1},M_{\omega_n}^{2},\dots ,M_{\omega_n}^{k-1},\dots ),\ n>1,
$$
where $M_{\omega_n}^{k-1}$ is at the $k$th position in the sequence.
\end{prop}

Indeed, let us assume that $M_{[n]}^{k-1}$ is at the $k$th
position in the sequence $(RX)^{[n]}$. Then, the $k$th component
of $(RX)^{[n]} \cdot X$ reads
$\sum_{0<i_1<\dots<i_n<k}x_{i_1}\dots x_{i_n}x_k$ and the $k$th
component of $R((RX)^{[n]}\cdot X)$ reads
$\sum_{i=1}^{k-1}\sum_{0<j_1<\dots<j_n<i}x_{j_1} \dots x_{j_n}
x_i=M_{[n+1]}^{k-1}$. The identity for $(RX)^{\{n\}}$ follows by
symmetry.

\begin{cor}
The elements $(RX)^{[n]}$ generate freely an associative
subalgebra of $A$. Similarly, the elements $(RX^i)^{[n]}$ generate
freely an associative subalgebra of $A^{(n)}$.
\end{cor}

Let us sketch the proof. Noncommutative monomials over $X$ such as
the ones appearing in the expansion of $M_{1,3,3,2}$ are naturally
ordered by the lexicographical ordering $<_L$, so that, for
example, $x_2x_6x_6x_4<_L x_2x_7x_7x_5$. Let us write, for any
noncommutative polynomial $P$ in $T(X)$, $Sup(P)$ for the highest
noncommutative monomial for the lexicographical ordering appearing
in the expansion of $P$, so that, for example, for $k \geq n$,
$Sup(M_{[n]}^{k}) =  x_{k-n+1} \dots x_k$, or
$Sup(x_2x_6x_7+x_2x_7x_5) = x_2x_7x_5$. For two such polynomials
$P$ and $Q$, we write $P<_LQ$ when $Sup(P)<_LSup(Q)$.

Let us consider now a noncommutative polynomial $Q$ in the
$(RX)^{[n]}$ with non trivial coefficients, and let us prove that
it is not equal to $0$ in $A$. For degree reasons, we may first
assume that $Q$ is homogeneous, that is, that $Q$ can be written
$$
    Q=\sum\limits_k
      \sum\limits_{n_1 + \dots + n_k=p}
      \alpha_{n_1,\dots ,n_k} (RX)^{[n_1]} \cdot \dots \cdot (RX)^{[n_k]},
$$
Then, the corollary follows from the observation that, for $l \gg
p$, $Sup(M_{[n_1]}^{l} \cdot \dots \cdot M_{[n_k]}^{l})
> Sup(M_{[m_1]}^{l}\cdot \dots \cdot M_{[m_j]}^{l})$ with
$n_1 + \dots + n_k = m_1 + \dots + m_j$ if and only if the
sequence $(n_1, \dots ,n_k)$ is less than the sequence $(m_1,
\dots ,m_j)$ in the lexicographical ordering. Indeed, let us
assume that the two sequences are distinct, and let $j$ be the
lowest index such that $m_j \not= n_j$, then $Sup(M_{[n_j]}^{l}) =
x_{l-n_j+1}\dots x_{l}$ whereas $Sup(M_{[m_j]}^{l})= x_{l-m_j+1}
\dots x_{l}$, so that, in particular $l-n_j+1>l-m_j+1$ if and only
if $m_j>n_j$. We can then conclude using the following obvious but
useful lemma.

\begin{lem}
For any $x,y$ homogeneous noncommutative polynomials in $T(X)$ and
$z,t$ in $T(X)$, due to the properties of the lexicographical
ordering, we have:
$$
    x <_L y \Rightarrow x \cdot z <_L y \cdot z
    \ \quad\ \
    {\rm{and}}
    \ \quad\ \
    z <_L t \Rightarrow x \cdot z <_L x \cdot t.
$$
\end{lem}

The proof goes over to $A^{(n)}$, provided one chooses a suitable
ordering on the elements of $X^1\coprod \dots \coprod X^n$, for
example the order extending the order on the $X^i$ and such that
$x_m^i<x_n^j$ whenever $i<j$.

\begin{cor}
The elements $(RX)^{[n]}$ generate freely an associative
subalgebra of $A$ for the double Rota--Baxter product
$\ast$.
\end{cor}
Here, we abbreviate the notation for the weight one double product
$\ast_1$ to $\ast$. The same assertion (and its proof) holds
{\textit{mutatis mutandis}} for the $(RX^i)^{[n]}$ and $A^{(n)}$.

The corollary follows from the previous lemma and the observation
that, for any sequence $(n_1,\dots ,n_k)$, $Sup(\{(RX)^{[n_1]}
\cdot \dots \cdot (RX)^{[n_k]}\}_{l}) = Sup(\{(RX)^{[n_1]} \ast
\dots \ast (RX)^{[n_k]} \}_l)$, where the lower $l$ indicates the
order of the component in the sequences.

Let us, once again, sketch the proof that relies on the usual
properties of the lexicographical ordering and the definition of
$(RX)^{[n]}$. Notice first that, for $k\gg n$,
$Sup(R((RX)^{[n]})_k) < Sup((RX)^{[n]}_k)$. Indeed, $(RX)^{[n]}_k
= M_{[n]}^{k-1}$, whereas
$$
    R((RX)^{[n]})_k = (RX)^{[n]}_1 + \dots + (RX)^{[n]}_{k-1}
                    = M_{[n]}^{1} + \dots + M_{[n]}^{k-2}
$$
so that
 \allowdisplaybreaks{
\begin{eqnarray*}
    Sup\bigl(R((RX)^{[n]})_k\bigr) &=& Sup\bigl(M_{[n]}^{k-2}\bigr)\\
                         &=& x_{k-n-1} \dots x_{k-2} <_L x_{k-n} \dots x_{k-1}
                          = Sup\bigl(M_{[n]}^{k-1}\bigr) = Sup\bigl(((RX)^{[n]})_k\bigr)
\end{eqnarray*}}

The same argument shows that, more generally, for any element
$Y=(y_1,\dots ,y_n,\dots )$ in $A$ satisfying the lexicographical
growth condition $Sup(y_i) <_L Sup(y_{i+1})$ for all $i \in
\N^\ast$, we have $Sup(y_{i-1}) = Sup(R(Y)_i) <_L Sup(Y_i) =
Sup(y_i)$ for all $i \in\N^\ast$. This property is therefore stable under
the map $R$ and is (up to neglecting the zero entries in the
sequences) common to all the elements we are going to consider. It
applies in particular to $(RX)^{[n_1]}\cdot \dots \cdot
(RX)^{[n_k]}$ and $(RX)^{[n_1]} \ast \dots \ast (RX)^{[n_k]}$; the
verification follows from the same line of arguments and is left
to the reader.

In the end, we have:
 \allowdisplaybreaks{
\begin{eqnarray*}
    (RX)^{[n_1]} \ast \dots \ast (RX)^{[n_k]}
    &=& (RX)^{[n_1]} \cdot \bigl((RX)^{[n_2]} \ast \dots \ast (RX)^{[n_k]}\bigr)\\
    & & + R\bigl((RX)^{[n_1]}\bigr) \cdot ((RX)^{[n_2]} \ast \dots
       \ast (RX)^{[n_k]})\\
    & &   + (RX)^{[n_1]} \cdot R\bigl((RX)^{[n_2]} \ast \dots \ast (RX)^{[n_k]}\bigr)
\end{eqnarray*}}
from which we deduce by recursion on $k$ and for $l \gg 0$:
 \allowdisplaybreaks{
\begin{eqnarray*}
    Sup\bigl( ((RX)^{[n_1]} \ast \dots \ast (RX)^{[n_k]})_l\bigr)
    &=& Sup\bigl(((RX)^{[n_1]} \cdot ((RX)^{[n_2]} \ast \dots \ast (RX)^{[n_k]}))_l\bigr)\\
    &=& Sup\bigl(((RX)^{[n_1]} \cdot \dots \cdot (RX)^{[n_k]})_l\bigr)
\end{eqnarray*}}

Let us conclude by proving that these constructions give rise to a
model for free Rota--Baxter algebras.

\begin{thm}
The RB subalgebra $\mathcal R$ of $A$ generated by $X$ is a free
RB algebra on one generator. More generally, the RB subalgebra
${\mathcal R}^{(n)}$ of $A^{(n)}$ generated by $X^1,\ldots,X^n$ is
a free RB algebra on $n$ generators.
\end{thm}

Our proof is inspired by the one in
Rota--Smith~\cite{rotasmith1972}, but the adaptation to the
noncommutative setting requires some care.

Let us write $\mathcal F$ for the free noncommutative RB algebra
on one generator $Y$, so that, by the universal properties of free
algebras, the map sending $Y$ to $X$ induces a surjective RB map
from $\mathcal F$ to $\mathcal R$ (recall that the latter is
generated by $X$).

Let us call $\End$-algebra any associative algebra $V$ provided
with a linear endomorphism $T_V \in \End(V)$, with the obvious
notion of morphisms of $\End$-algebras, so that a $\End$-algebra
morphism $f$ from $V$ to $W$ satisfies $f \circ T_V = T_W \circ
f$. Let us write $\mathcal L$ for the free $\End$-algebra on one
generator. We now write $Z$ for the generator: the elements of
$\mathcal L$ are linear combinations of all the symbols obtained
from $Z$ by iterative applications of the endomorphism $T$ and of
the associative product. The elements look like
$ZT^2(T(Z)T^3(Z))$, and so on. We write $M$ for the set of these
symbols and call them $\mathcal L$-monomials.

A RB algebra $B$ is a $\End$-algebra together with extra
(Rota--Baxter) relations on $T_B=R$. In particular, there is a
unique natural $\End$-algebra map from $\mathcal L$ to an
arbitrary RB algebra on one generator (mapping $Z$ to that
generator) and in particular a unique map to $\mathcal F$ and
$\mathcal R$ sending $Z$ to $Y$, resp. $Z$ to $X$. The map to
$\mathcal F$ factorizes the map to $\mathcal R$.

Proving that $\mathcal F$ and $\mathcal R$ are isomorphic as RB
algebras, that is, that $\mathcal R$ is a free RB algebra on one
generator amounts to prove that the kernel ---say $Ker(F)$--- of
the map from $\mathcal L$ to $\mathcal F$ is equal to the kernel
---say $Ker(U)$--- of the map from $\mathcal L$ to $\mathcal R$.

For any $l \in \mathcal L$, which can be written uniquely as a
linear combination of $\mathcal L$-monomials, we write $Max(l)$
for the maximal number of $T$s occurring in the monomials (with
the obvious conventions for the powers of $T$, so that for
example: $Max(ZT^2(ZT(Z))+Z^3T^2(Z)Z)=3$).

We say that an element $\alpha$ of $M$ is
{\it{elementary}} if and only if it can be written either $Z^i$,
$i \geq 0$ or as a product $Z^{i_1} \cdot T(b_1) \cdot
Z^{i_2}\cdot \dots \cdot T(b_k) \cdot Z^{i_{k+1}}$, where the
$b_i$'s are elementary and $i_2,\dots ,i_k$ are strictly positive
integers ($i_1$ and $i_{k+1}$ may be equal to zero); the
definition of elementariness makes sense by induction on $Max(
\alpha )$.

\begin{lem}
Every $l \in \mathcal L$ is of the form $l=r+s$ where $F(s)=0$ and
$r$ is a sum of elementary monomials.
\end{lem}

It is enough to prove the lemma for $l = t \in M$. If $t$
is not elementary, then $t$ has at least in its expansion a
product of two consecutive factors of the form $T(c) \cdot T(d)$.
However, since $\mathcal L$ is a Rota--Baxter algebra, the
relation
$$
    T(c \cdot T(d) + T(c) \cdot d + c \cdot d) - T(c)\cdot T(d) \in Ker(L)
$$
holds, and $t$ can be rewritten, up to an element in $Ker(L)$, by
substituting $T(c \cdot T(d)+T(c) \cdot d + c \cdot d)$ to
$T(c)\cdot T(d)$ in its expansion. Notice that $Max(c \cdot T(d)+T(c)
\cdot d + c \cdot d) < Max (T(c)\cdot T(d))$. The proof follows by a
joint induction on the number of such consecutive factors and on
$Max(T(c)\cdot T(d))$. In other words, products $T(c)\cdot T(d)$
can be iteratively cancelled from the expression of $t$ using the
RB fundamental relation.

Let  us show now that, with the notation of the lemma, $U(r)=0$
implies $F(r)=0$, from which, since we already know that $Ker(F)
\subset Ker(U)$, the freeness property will follow.

We actually claim the stronger property that, for $p$ large enough
and for $\mu \not= \mu'$ elementary, $Sup(U(\mu)_p) \not=
Sup(U(\mu')_p)$, from which the previous assertion will follow.
Indeed, if $\mu
 = (Z^{i_1} \cdot T(b_1) \cdot Z^{i_2} \cdot \dots \cdot T(b_k) \cdot
Z^{i_{k+1}})$ with the $b_i$ elementary, then
 \allowdisplaybreaks{
\begin{eqnarray*}
   \lefteqn{Sup(U(\mu )_p)=x_p^{i_1}(Sup(U(T(b_1)))_{p})x_p^{i_2}\dots x_p^{i_k}(Sup(U(T(b_k)))_{p})
   x_p^{i_{k+1}}}\\
    &=&x_p^{i_1}(Sup(R(U(b_1)))_{p})x_p^{i_2}\dots x_p^{i_k}(Sup(R(U(b_k)))_{p}) x_p^{i_{k+1}}
\end{eqnarray*}}

the last identity follows  since $U$ is an $\End$-algebra map. Since
$U(\tau)$ for $\tau$ elementary satisfies the lexicographical
growth condition (as may be checked by induction), we have
$Sup(R(U(\tau))_i)=Sup(U(\tau)_{i-1})$, so that
$$
    Sup(U(\mu )_p) = x_p^{i_1} (Sup(U(b_1)))_{p-1}) x_p^{i_2} \dots x_p^{i_k}(Sup((U(b_k))_{p-1}) x_p^{i_{k+1}}.
$$
The proof follows by induction on $Max(\mu)$.

The proof goes over to an arbitrary number of generators, provided
one defines the suitable notion of elementary monomials in the
free $\End$-algebra ${\mathcal L}^{(n)}$ on $n$ generators
$Z_1,\dots ,Z_n$: these are the elements of ${\mathcal L}^{(n)}$
that can be written either as noncommutative monomials in the
$Z_i$, or as a product $a_1T(b_1)a_2\dots T(b_n)a_{n+1}$, where
the $b_i$ are elementary and the $a_i$ noncommutative monomials in
the $Z_i$ (nontrivial whenever $1<i<n+1$).

Notice the following interesting corollary of our previous
computations.

\begin{cor}
The images of the elementary monomials of $\mathcal L$ in
$\mathcal R$ form a basis (as a vector space) of the free
Rota--Baxter algebra on one generator.
\end{cor}

The same assertion holds for the free Rota--Baxter algebra on $n$
generators.


\section{Rota--Baxter algebras and NCSF}
\label{sect: RBandNCSF}

In the present section, we associate to the free RB algebra on one
generator a Hopf algebra naturally isomorphic to the Hopf algebra
of noncommutative symmetric functions or, equivalently, to the
descent algebra. The reasons for the introduction of this Hopf
algebra will become clear in the next section. Let us simply
mention at this stage that these notions will provide the right
framework to extend to noncommutative RB algebras and
noncommutative symmetric functions the classical results of Rota
and Smith relating commutative RB algebras and symmetric
functions~\cite{rotasmith1972}.

Recall that the algebra $NCQSym$ of quasi-symmetric functions in
noncommutative variables introduced in the previous section is
naturally provided with a Hopf algebra
structure~\cite{bergeron2005}. On the elementary quasi-symmetric
functions $M_{[n]}$, the coproduct $\Delta$ acts as on a sequence
of divided powers
$$
    \Delta \bigl(M_{[n]}\bigr)=\sum\limits_{i=0}^n M_{[i]}\otimes M_{[n-i]}.
$$
The same argument as in the previous section shows that the
$M_{[n]}$ generate a free subalgebra of $NCQSym$. In the end, the
$M_{[n]}$s form a sequence of divided powers in a free associative
sub-algebra of $NCQSym$, and this algebra is isomorphic to the
algebra of NCSF (which is, by its very definition a Hopf algebra
freely generated as an associative algebra by a sequence of
divided powers~\cite{gelfand1995}, and therefore is naturally
isomorphic to the descent algebra: see Prop.~\ref{propfree} and
the description of the Hopf algebra structure on the descent
algebra in the same section).

The same construction of a Hopf algebra structure goes over to the
algebras introduced in the previous section, that is, to the free
algebras over the $(RX)^{[n]}$ for the $\cdot$ and $\ast$
products. As a free algebra over the $(RX)^{[n]}$, the first
algebra is naturally provided with a cocommutative Hopf algebra
structure for which the $(RX)^{[n]}$s form a sequence of divided
powers, that is
\begin{equation}
\label{SpitzerCoproduct1}
    \Delta \bigl((RX)^{[n]}\bigr)
    =\sum\limits_{0 \leq m \leq n} (RX)^{[m]} \otimes (RX)^{[n-m]}.
\end{equation}
This is the structure inherited from the Hopf algebra structure on
$NCQSym$. We will be particularly interested in this Hopf algebra,
that is the algebra freely generated by the $(RX)^{[n]}$ for the
$\cdot$ product, viewed as a subalgebra of $A$ and as a Hopf
algebra. We call it the {\textit{free noncommutative Spitzer
(Hopf) algebra on one generator}} or, for short, the Spitzer
algebra, and write it $\mathcal{S}$.

When dealing with the Rota--Baxter double product, $\ast$, the
right subalgebra to consider, as will appear below, is not the
free algebra generated by the $(RX)^{[n]}$ but the free algebra
freely generated by the $(RX)^{[n]} \cdot X$. We may also consider
this algebra as a Hopf algebra by requiring the free generators to
form a sequence of divided powers, that is by defining the
coproduct by
\begin{equation}
\label{SpitzerCoproduct2}
    {\Delta}_\ast \bigl( (RX)^{[n-1]} \cdot X \bigr)
    = (RX)^{[n-1]}\cdot X \otimes 1
     + \sum\limits_{m \leq n-2} (RX)^{[m]} \cdot X \otimes (RX)^{[n-m-2]} \cdot X
     + 1 \otimes (RX)^{[n-1]}\cdot X
\end{equation}
We will also investigate briefly this second structure, strongly
related with the Hopf algebra structure on free dendriform
dialgebras of~\cite{chapoton2000,ronco2000}. We call it the
{\textit{double Spitzer algebra}} ---and write it $\mathcal{C}$.

To understand first the structure of the Spitzer algebra,
$\mathcal{S}$, recall that we write $\tilde R$ for $- \theta \id -
R$ and $(\tilde{R}a)^{\{n\}}$ respectively $(\tilde{R}a)^{[n]}$
for the corresponding iterated operator.

\begin{lem}
\label{cor:RBantipode} The action of the antipode $S$ in the
Spitzer algebra, $\mathcal{S}$, is given by
$$
    S\bigl((RX)^{[n]}\bigr) = -R\bigl(X \cdot ({\tilde R}X)^{\{n-1\}}\bigr).
$$
\end{lem}

Indeed, the Spitzer algebra is naturally a graded Hopf algebra.
The series $F:=\sum_{n \ge 0} (RX)^{[n]}$ is a group-like element
in the Hopf algebra. The inverse series follows from Atkinson's
formula~\ref{inverseAtkinson}, and gives the action of the
antipode on the terms of the series. Since
$$
    F^{-1} = 1 - R\bigl( X \cdot (\sum\limits_{n \ge 0}({\tilde R}X)^{\{n\}})\bigr),
$$
the corollary follows.

\begin{cor}
\label{cor:RBdoubleAntipode} The action of the antipode $S$ in the
double Spitzer algebra, $\mathcal{C}$, is given by
$$
    S\bigl((RX)^{[n]}\cdot X\bigr) = - \bigl(X \cdot ({\tilde R}X)^{\{n\}}\bigr).
$$
\end{cor}

The proof will illustrate the links between the two Hopf algebras,
$\mathcal S$ and $\mathcal C$. Recall that, on any RB algebra, we
have, by the very definition of the $\ast_\theta$ product
$$
    R(x) \cdot R(y) = R(x \ast_\theta y)
    \quad {\rm{ and }}\quad
    \tilde{R}(x) \cdot \tilde{R}(y) = -\tilde{R}(x \ast_\theta y).
$$
Recall also that, in the algebra of series $A$ (in which the
Spitzer algebra and the double Spitzer algebra can be embedded),
the operator $R$ can be inverted on the left ---that is, if
$Y=(0,y_1,y_2,\dots,y_n,\dots)=R(U)$, then
$U=(y_1,y_2-y_1,y_3-y_2,\dots)$.

It follows from this observation that the RB operator $R$ induces
an isomorphism of free graded algebras between the double Spitzer
algebra and the Spitzer algebra. That is, for any sequence of
integers $i_1,\dots,i_k$, we have
$$
    R\Bigl(\bigl((RX)^{[i_1]} \cdot X\bigr) \ast \bigl((RX)^{[i_2]} \cdot X\bigr) \ast
     \cdots \ast \bigl((RX)^{[i_k]} \cdot X\bigr)\Bigr)
    =(RX)^{[i_1+1]} \cdot (RX)^{[i_2 +1]} \cdot \cdots \cdot (RX)^{[i_k+1]}
$$
The isomorphism is extended by the identity to the scalar (that is
to the zero degree components of the two Hopf algebras).

Since this isomorphism maps the generators of $\mathcal{C}$ to the
generators of $\mathcal S$, and since both families of generators
form a sequence of divided powers in their respective Hopf
algebras, we obtain that, writing again $S$ for the antipode in
$\mathcal C$
 \allowdisplaybreaks{
\begin{eqnarray*}
    R\Bigl(S\bigl(\sum\limits_{n > 0} (RX)^{[n-1]}\cdot X\bigr)\Bigr)
    &=& R\Bigl(\bigl(\sum\limits_{n > 0} (RX)^{[n-1]}\cdot X\bigr)^{-1}\Bigr)\\
    &=&\bigl(\sum\limits_{n > 0}(RX)^{[n]}\bigr)^{-1}
     =\sum \limits_{n > 0} -R\bigl(X\cdot ({\tilde R}X)^{\{n-1\}}\bigr)
\end{eqnarray*}}
so that, eventually
$$
    S\bigl((RX)^{[n]} \cdot X\bigr) = - X \cdot ({\tilde R}X)^{\{n\}}.
$$

\begin{cor}
The free $\ast$ subalgebras of $A$ generated by the
$(RX)^{[n]}\cdot X$ and by the $X \cdot ({\tilde R}X)^{\{n\}}$
identify canonically. The antipode exchanges the two families of
generators. In particular, the $X \cdot ({\tilde R}X)^{\{n\}}$
also form a sequence of divided powers in the double Spitzer
algebra.
\end{cor}


\section{The Bohnenblust--Spitzer formula and the Dynkin idempotent}
\label{sect:BSandDynkin}

As already alluded at, one surprising conclusion of the present
article is that {\textit{the same machinery}} that one uses to
derive fundamental identities in the theory of free Lie algebras
can be used to recover the already known formulas for commutative
Rota--Baxter algebras but, moreover, can be used to prove new
formulas in the noncommutative setting. These results rely on the
computation of the action of the Dynkin operator on the generators
of the Spitzer and of the double Spitzer Hopf algebras.

Let us now introduce the definition of the iterated Rota--Baxter
left and right pre-Lie brackets in a RB algebra $(B,R)$ of weight
$\theta$.
 \allowdisplaybreaks{
\begin{eqnarray}
    \mathfrak{l}_{\theta}^{(n)}(a_1,\dots,a_n) &:=&
    \Bigl( \cdots \bigl( (a_1 \triangleright_\theta a_2) \triangleright_\theta a_3 \bigr)
    \cdots \triangleright_\theta a_{n-1} \Bigr) \triangleright_\theta a_n
\label{leftRBpreLie}\\
    \mathfrak{r}_{\theta}^{(n)}(a_1,\dots,a_n) &:=&
    a_1 \triangleleft_\theta \Bigl( a_2 \triangleleft_\theta \bigl( a_3 \triangleleft_\theta
    \cdots (a_{n-1} \triangleleft_\theta a_n) \bigr) \cdots \Bigr)
\label{rightRBpreLie}
\end{eqnarray}}
for $n>0$ and
$\mathfrak{l}_{\theta}^{(1)}(a):=a=:\mathfrak{r}_{\theta}^{(1)}(a)$.
For fixed $a \in B$ we can write compactly for $n>0$
 \allowdisplaybreaks{
\begin{eqnarray}
    \mathfrak{l}_{\theta}^{(n+1)}(a) = \bigl(\mathfrak{l}_{\theta}^{(n)}(a)\bigr)\triangleright_\theta a
    \quad\ {\rm{ and }} \quad\
    \mathfrak{r}_{\theta}^{(n+1)}(a) = a \triangleleft_\theta \bigl(\mathfrak{r}_{\theta}^{(n)}(a)\bigr).
\label{def:pre-LieWords}
\end{eqnarray}}
We call those expressions left respectively right RB pre-Lie
words. Let us now define
 \allowdisplaybreaks{
\begin{eqnarray}
    \mathfrak{L}_{\theta}^{(n+1)}(a):=R\bigl(\mathfrak{l}_{\theta}^{(n+1)}(a)\bigr)
    \quad\ {\rm{ and }} \quad\
    \mathfrak{R}_{\theta}^{(n+1)}(a):=R\bigl(\mathfrak{r}_{\theta}^{(n+1)}(a)\bigr).
\label{def:LamCtheta}
\end{eqnarray}}
For $B$ commutative $\mathfrak{L}_{\theta}^{(n)}(a)=
(-\theta)^{n-1}R(a^n)$ and $\mathfrak{R}_{\theta}^{(n)}(a)=
(-\theta)^{n-1}R(a^n)$. For $(B,R)$ being of weight $\theta = 0$
the left (right) pre-Lie product~(\ref{lem:leftRBpreLie}) reduces
to $a \triangleright_0 b = ad_{R(a)}(b)$ (and $a \triangleleft_0 b
= -ad_{R(b)}(a)$), so that
 \allowdisplaybreaks{
\begin{eqnarray}
  \mathfrak{L}_0^{(n+1)}(a) = R\Bigl( \bigl[R\bigl( \cdots \bigl[ R([ R(a),a ]),a \bigr] \cdots \bigr),a \bigr]\Bigr)
                            = - R\Bigl( ad_{a}(\mathfrak{L}_0^{(n)}(a)) \Bigr)
\label{def:LamsC0}
\end{eqnarray}}
for $n>0$ and analogously for the right RB pre-Lie words. Now, in
the context of the weight $\theta=1$ Rota--Baxter algebra $(A,R)$
we find the following proposition.

\begin{prop}
\label{prop:RBDynkin} The action of the Dynkin operator on the
generators $(RX)^{[n]}$ of the Spitzer algebra $\mathcal{S}$
(respectively on the generators of the double Spitzer algebra
$\mathcal{C}$) is given by:
$$
    D\bigl((RX)^{[n]}\bigr) = \mathfrak{L}_{1}^{(n)}(X)
                            = R(\mathfrak{l}_{1}^{(n)}(X))
$$
(respectively by $D((RX)^{[n]} \cdot X) =
\mathfrak{l}_{1}^{(n+1)}(X)$).
\end{prop}

Due to the existence of the Hopf algebra isomorphism induced by
the map $R$ between the Spitzer algebra and the double Spitzer
algebra, the two assertions are equivalent. Let us prove the
proposition by induction on $n$ for the double Spitzer algebra. We
denote by~$\pi_*$ the product on~$\mathcal C$. Using $Y(1)=0$ we
find for $n=0$
$$
    D(X) = (S \star Y)(X)=\pi_* \circ (S \otimes Y)(X \otimes 1 + 1 \otimes X)
                         = X
                         = \mathfrak{l}_{1}^{(1)}(X).
$$
Recall that $Y(X)=X$ and $Y((RX)^{[n-1]}X)= n \ (RX)^{[n-1]}X=
Y((RX)^{[n-1]}) \cdot X + (RX)^{[n-1]}\cdot X$. Let us also
introduce a useful notation and write $(RX)^{[n-1]} \cdot X =:
w^{(n)}$, so that $ w^{(n+1)} = R(w^{(n)}) \cdot X$, $w^{(1)}=1$
and $w^{(0)}=1$. We obtain
 \allowdisplaybreaks{
\begin{eqnarray*}
 \lefteqn{D(w^{(n)}) = D\bigl((RX)^{[n-1]} \cdot X\bigr)
                        = (S \star Y)\bigl(w^{(n)}\bigr)}\\
        &=& \pi_* \circ (S \otimes Y)
            \sum\limits_{p =0}^{n} w^{(p)} \otimes w^{(n-p)}
         = \sum\limits_{p =0}^{n-1} S\bigl(w^{(p)}\bigr) \ast
                                     Y\bigl( R(w^{(n-1-p)}) \cdot X \bigr)\\
        &=& \sum\limits_{p =0}^{n-1} S\bigl(w^{(p)}\bigr) \ast
                                     \bigl(Y\bigl( R(w^{(n-1-p)})\bigr) \cdot X\bigr)
          + \sum\limits_{p =0}^{n-1} S\bigl(w^{(p)}\bigr) \ast
                                      \bigl(R(w^{(n-1-p)}) \cdot X\bigr)\\
        &=& \sum\limits_{p =0}^{n-1} S\bigl(w^{(p)}\bigr) \ast
                                     \bigl(Y\bigl( R(w^{(n-1-p)})\bigr) \cdot X\bigr)
          + \sum\limits_{p =0}^{n} S\bigl(w^{(p)}\bigr) \ast
                                      w^{(n-p)} -  S\bigl(w^{(n)}\bigr)\\
        &=& \sum\limits_{p =1}^{n-1} S\bigl(w^{(p)}\bigr) \ast
                                     \bigl(Y\bigl( R(w^{(n-1-p)})\bigr) \cdot X\bigr)
          + Y\bigl( R(w^{(n-1)})\bigr) \cdot X - S\bigl(w^{(n)}\bigr)
\end{eqnarray*}}
where, since the antipode $S$ is the convolution inverse of the
identity, the term $(S \star \id)(w^{(n)})$ on the right hand side
cancels. Hence, using the general RB identity
$$
    a \ast \bigl(R(b) c\bigr) = R(a) R(b) c - a{\tilde R}\bigl(R(b)c\bigr)
                           = R\bigl(a \ast b\bigr)c - a{\tilde R}\bigl(R(b)c\bigr)
$$
we find immediately
 \allowdisplaybreaks{
\begin{eqnarray*}
    D\bigl(w^{(n)}\bigr)
    &=&  \sum\limits_{p =1}^{n-1} R\Bigl( S\bigl(w^{(p)}\bigr) \ast
                                     Y(w^{(n-1-p)})\Bigr) \cdot X
       - \sum\limits_{p =1}^{n-1} S\bigl(w^{(p)}\bigr) \cdot
                                  \tilde{R}\Bigl( R\bigl(Y(w^{(n-1-p)})\bigr) \cdot X \Bigr)\\
    & & + Y\bigl( R(w^{(n-1)})\bigr) \cdot X - S\bigl(w^{(n)}\bigr)\\
    &=&  \sum\limits_{p =0}^{n-1} R\Bigl( S\bigl(w^{(p)}\bigr) \ast
                                     Y(w^{(n-1-p)})\Bigr) \cdot X
       - \sum\limits_{p =1}^{n-1} S\bigl(w^{(p)}\bigr) \cdot
                                  \tilde{R}\Bigl( R\bigl(Y(w^{(n-1-p)})\bigr) \cdot X \Bigr)
       - S\bigl(w^{(n)}\bigr)\\
    &=& R\bigl( (S \star Y) (w^{(n-1)})\bigr) \cdot X -
        S\bigl(w^{(n)}\bigr)
       - \sum\limits_{p =1}^{n-1} X \cdot \tilde{R}\Bigl(S\bigl(w^{(p-1)}\bigr)\Bigr)
                                    \cdot \tilde{R}\Bigl( R\bigl(Y(w^{(n-1-p)})\bigr) \cdot X \Bigr)
\end{eqnarray*}}
where we used that
$$
   S\bigl(w^{(p)}\bigr) = S\bigl((RX)^{[p-1]} \cdot X \bigr)
                        = - X \cdot ({\tilde R}X)^{\{p-1\}}
                        = - X \cdot {\tilde R}\bigl( X \cdot (\tilde{R}X)^{\{p-2\}}\bigr)
                        =   X \cdot {\tilde R}\bigl( S(w^{(p-1)}) \bigr).
$$
Recall that $\tilde R$ is a Rota--Baxter operator as well, such
that $\tilde{R}(a \ast_{\theta} b) = -\tilde{R}(a)\tilde{R}(b)$.
This leads to
 \allowdisplaybreaks{
\begin{eqnarray*}
    D\bigl(w^{(n)}\bigr)
    &=& R\bigl( D(w^{(n-1)})\bigr) \cdot X -
        S\bigl(w^{(n)}\bigr)
       + \sum\limits_{p =1}^{n-1} X \cdot \tilde{R}\Bigl(
                                     S\bigl(w^{(p-1)}\bigr)
                                      \ast
                                     \bigl(Y\bigl(R(w^{(n-1-p)})\bigr) \cdot X \bigr)\Bigr)\\
    &=& R\bigl( \mathfrak{l}_{1}^{(n-1)})\bigr) \cdot X
       - X \cdot \tilde{R}\bigl( S\bigl(w^{(n-1)}\bigr)
       + \sum\limits_{p =1}^{n-1} X \cdot \tilde{R}\Bigl(
                                     S\bigl(w^{(p-1)}\bigr)
                                      \ast
                                     \bigl(Y\bigl(R(w^{(n-1-p)}) \cdot X \bigr)\Bigr)\\
    & &   - \sum\limits_{p =1}^{n-1} X \cdot \tilde{R}\Bigl(
                                     S\bigl(w^{(p-1)}\bigr)
                                      \ast
                                     \bigl(R(w^{(n-1-p)}) \cdot Y(X) \Bigr)\\
    &=& R\bigl( \mathfrak{l}_{1}^{(n-1)}(X)\bigr) \cdot X
       + \sum\limits_{p =0}^{n-2} X \cdot \tilde{R}\Bigl(
                                     S\bigl(w^{(p)}\bigr)
                                      \ast
                                     Y\bigl(w^{(n-1-p)}\bigr)\Bigr)
                                     - \sum\limits_{p =0}^{n-1} X \cdot \tilde{R}\Bigl(
                                     S\bigl(w^{(p)}\bigr)
                                      \ast
                                     w^{(n-1-p)}\Bigr)\\
    &=& R\bigl( \mathfrak{l}_{1}^{(n-1)}(X)\bigr) \cdot X
       + \sum\limits_{p =0}^{n-1} X \cdot \tilde{R}\Bigl(
                                     S\bigl(w^{(p)}\bigr)
                                      \ast
                                     Y\bigl(w^{(n-1-p)}\bigr)\Bigr)
            - X \cdot \tilde{R}\bigl( (S \star \id) (w^{(n-1)})\bigr)\\
    &=& R\bigl( \mathfrak{l}_{1}^{(n-1)})\bigr) \cdot X
       + X \cdot \tilde{R}\bigl((S \star Y)(w^{(n-1)})\bigr)\\
    &=& R\bigl( \mathfrak{l}_{1}^{(n-1)}(X)\bigr) \cdot X +  X \cdot \tilde{R}\bigl( D(w^{(n-1)})\bigr)\\
    &=& R\bigl( \mathfrak{l}_{1}^{(n-1)}(X)\bigr) \cdot X +  X \cdot \tilde{R}\bigl( \mathfrak{l}_{1}^{(n-1)}(X)\bigr)
     =\mathfrak{l}_{1}^{(n-1)}(X) \triangleright_1 X = \mathfrak{l}_{1}^{(n)}(X)
\end{eqnarray*}}
We used once again that the antipode $S$ is the convolution
inverse of the identity, implying that $- X \cdot \tilde{R}\bigl(
(S \star \id) (w^{(n-1)})\bigr)=0$. All this immediately implies
the following important theorem.

\begin{thm} \label{thm:ncBS1}
We have the following identity in the Spitzer algebra
$\mathcal{S}$:
\begin{equation}
\label{eq:key1}
    (RX)^{[n]}=\sum\limits_{i_1+\cdots+i_k = n, \atop  i_1,\dots,i_k>0}
    \frac{\mathfrak{L}_{1}^{(i_1)}(X)\cdot \ \cdots\ \cdot \mathfrak{L}_{1}^{(i_k)}(X)}{i_1(i_1+i_2)\cdots(i_1+\cdots+i_k)}
\end{equation}
\end{thm}

The theorems follows readily from Theorem~\ref{thm:invDynkinRB} by
applying the formula for the inverse of the Dynkin operator in
Proposition~\ref{prop:RBDynkin}. We obtain the equivalent
expansion in the double Spitzer algebra $\mathcal{C}$

\begin{cor}
We have, in the double Spitzer algebra $\mathcal{C}$:
$$
    (RX)^{[n-1]} \cdot X=
    \sum\limits_{i_1+ \cdots +i_k=n, \atop  i_1,\dots,i_k>0}
    \frac{\mathfrak{l}_{1}^{(i_1)}(X) \ast \cdots \ast
    \mathfrak{l}_{1}^{(i_k)}(X)}{i_1(i_1+i_2)\cdots(i_1+\cdots+i_k)}.
$$
\end{cor}

The reader should have no problem to verify the following
statements.

\begin{cor}
We have:
 \allowdisplaybreaks{
\begin{eqnarray*}
    (RX)^{\{n\}} &=&
    \sum\limits_{i_1+\cdots+i_k = n, \atop  i_1,\dots,i_k>0}
    \frac{\mathfrak{R}_{1}^{(i_k)}(X)\cdot \ \cdots\ \cdot
    \mathfrak{R}_{1}^{(i_1)}(X)}{i_1(i_1+i_2)\ldots(i_1+\cdots+i_k)}\\
     X \cdot (RX)^{\{n-1\}}&=&
    \sum\limits_{i_1+ \cdots +i_k=n, \atop  i_1,\dots,i_k>0}
    \frac{\mathfrak{r}_{1}^{(i_k)}(X) \ast \cdots \ast
    \mathfrak{r}_{1}^{(i_1)}(X)}{i_1(i_1+i_2)\cdots(i_1+\cdots+i_k)}.
\end{eqnarray*}}
with $\mathfrak{r}_{1}^{(i)}(X) $ and $\mathfrak{R}_{1}^{(i)}(X)$
defined in~(\ref{rightRBpreLie}) and (\ref{def:LamCtheta}),
respectively.
\end{cor}

At this point we may assume the Rota--Baxter algebra to be of
weight $\theta$, i.e. we replace the left-to-right bracketed
Rota--Baxter pre-Lie words $R(\mathfrak{l}_{1}^{(i_k)}(X))$ by
$R(\mathfrak{l}_{\theta}^{(i_k)}(X))$. C.S.~Lam discovered
in~\cite{lam1998}, see also~\cite{oteoros2000}, the weight zero
case of identity~(\ref{eq:key1}), that is, for
$\mathfrak{L}_0^{(i_k)}$.

Let us continue to follow closely Rota--Smith's
work~\cite{rotasmith1972} implying the natural extension to the
free Rota--Baxter algebra $\mathcal{R}$ in $n$ generators, i.e.,
sequences $X_1,\dots,X_n$. Now, working in the power series ring
$\mathcal{R}[t_1, \ldots, t_n]$ with $n$ commuting parameters
$t_1, \dots, t_n$, and replacing $X$ by $X_1t_1 + \dots + X_nt_n
\in \mathcal{R}[t_1, \ldots, t_n]$ in identity~(\ref{eq:key1}) of
Theorem~\ref{thm:ncBS1} we obtain a noncommutative generalization
of the classical Bohnenblust--Spitzer identity by comparing the
coefficients of the monomial $t_1 \dots t_n$ on both sides. We
arrive at the following identity for arbitrary RB algebras:

\begin{thm}\label{NCSpitzer}
Let $R$ be a Rota--Baxter operator on a Rota--Baxter algebra $A$
and $x_1,\ldots,x_n \in A$. Then:
 \allowdisplaybreaks{
\begin{eqnarray*}
    \sum_{\sigma\in S_n}
        R\Bigl(R\bigl(\cdots R(x_{\sigma_1}) x_{\sigma_2} \dots \bigr) x_{\sigma_n}\Bigr)
        =
        \sum\limits_{\pi \in \mathcal{OP}_n} \omega(\pi)
        \mathtt{L}_{\theta}(\pi_1)\cdot \ \cdots\  \cdot \mathtt{L}_{\theta}(\pi_k)
\end{eqnarray*}}
\end{thm}

The sum on the left hand side is over all permutations in $S_n$.
The sum on the right hand side is over all {\it{ordered}}
partitions, $\pi=[\pi_1]\dots[\pi_k]$, that is, sequences of its
disjoint subsets whose union is $[n]$. We denote by
$\mathcal{OP}_n=\sum_{i=1}^{n}\mathcal{OP}^k_n$ the set of all
ordered partitions and by $\mathcal{OP}^k_n$ the set of ordered
partitions of $[n]$ with $k$ blocks. We denote by $m_i:=|\pi_i|$
the number of elements in the block $\pi_i$ of partition $\pi$.
The coefficient function $\omega(\pi)$ is simply defined to be
$$
    \omega(\pi):=\bigl(m_1(m_1 + m_2)\cdots(m_1 + \cdots + m_k)\bigr)^{-1}.
$$
Finally, we define $\mathtt{L}_{\theta}(\pi_i)$, $[\pi_i]=[j_1
\ldots j_{m_i}]$, using the left-to-right bracketed RB pre-Lie
words of weight $\theta$ (\ref{leftRBpreLie}) in
$(\mathcal{R},R)$ by:
\begin{equation}
    \mathtt{L}_{\theta}(\pi_i) := \sum_{\sigma \in S_{m_i}}
             R\bigl({\mathfrak{l}}_{\theta}^{(m_i)}(X_{j_{\sigma_1}},\dots,X_{j_{\sigma_{m_i}}})\bigr).
\end{equation}
We recover identity~(\ref{eq:key1}), that is, for
$X=X_1=\ldots=X_n$, from the fact that the number of ordered
partitions of type $m_1 + \ldots + m_k=n$ is given by the
multinomial coefficient $n!(m_1! \dots m_k!)^{-1}$ and from the
fact that in that case we have
$$
      \mathtt{L}_{\theta}(\pi_i)=m_i!\mathfrak{L}_{\theta}^{(m_i)}(X).
$$

Let us turn now to the classical, commutative case and show how
the classical Spitzer identity can be recovered from the
noncommutative one. Recall first that, in the commutative case,
$a\triangleright_\theta b=-\theta ab$, so that, for $\pi_i$ as
above, we have
$$
      \mathtt{L}_{\theta}(\pi_i)=m_i!(-\theta)^{m_i-1}
      R\bigl(\prod_{j \in \pi_i}X_{j}\bigr)
$$
and we get
$$
    \sum\limits_{\sigma\in S_n}
    R(R(\ldots R(X_{\sigma_1})X_{\sigma_2}\dots )X_{\sigma_n})
    =\sum\limits_{\pi \in \mathcal{OP}_n}\frac{(m_1)!\cdots (m_k)!}
    {m_1 (m_1+m_2) \cdots (m_1+ \dots +m_k)}(-\theta )^{n-k} \prod\limits_{i=1}^k R(\prod\limits_{j\in\pi_i}X_j).
$$

\begin{lem}
We have, for any sequence $(m_1,\ldots,m_k)$:
$$
    \sum\limits_{\sigma\in S_k}\frac{1}{m_{\sigma_1} (m_{\sigma_1}+m_{\sigma_2}) \cdots (m_{\sigma_1}+ \cdots + m_{\sigma_k})}
    =\prod\limits_{i=1}^k\frac{1}{m_i}
$$
\end{lem}

Indeed, let us consider the integral expression
$$
    \int\limits_0^1 x_k^{m_k-1}dx_k \dots \int\limits_0^1 x_1^{m_1-1}dx_1 = \prod\limits_{i=1}^k\frac{1}{m_i}
$$
Recall that
$$
    \int\limits_0^1 x^pdx \int\limits_0^1 y^qdy = \int_0^1 x^p \int_0^x y^q dydx + \int_0^1 y^q \int_0^y x^p dxdy.
$$
The formula (a weight zero RB relation for the Riemann integral map) follows from the geometric decomposition of the square
into two triangles and generalizes to higher products of
integrals. In the general case, the hypercube in dimension $n$ is
divided into $n!$ simplices. We get:
 \allowdisplaybreaks{
\begin{eqnarray*}
    \int\limits_0^1 x_k^{m_k-1} dx_k \dots \int\limits_0^1 x_1^{m_1-1} dx_1
    &=&\sum\limits_{\sigma\in S_k} \int_0^1 x_{\sigma_1}^{m_{\sigma_1}-1}\int_0^{x_{\sigma_1}}
    x_{\sigma_2}^{m_{\sigma_2}-1} \dots \int_0^{x_{\sigma_{k-1}}}
    x_{\sigma_k}^{m_{\sigma_k}-1}dx_{\sigma_{k-1}}\ldots dx_{\sigma_1}\\
    &=&\sum\limits_{\sigma\in S_k} \frac{1}{m_{\sigma_1} (m_{\sigma_{1}}+m_{\sigma_{2}})  \cdots (m_{\sigma_1} + \cdots +m_{\sigma_k})}
\end{eqnarray*}}
which gives the expected formula.

This leads to the classical Bohnenblust--Spitzer
formula~\cite{rotasmith1972} of weight~$\theta$
\begin{equation}
    \sum_{\sigma\in S_n}
    R\Bigl( R\bigl( \cdots R(X_{\sigma_1})X_{\sigma_2} \cdots  \bigr) X_{\sigma_n} \Bigr)
    = \sum_{\pi \in \mathcal{P}_n} (-\theta)^{n -|\pi|}
      \prod_{\pi_i \in \pi}(m_i-1)! \ R\Bigl(\prod_{j \in \pi_i}X_j\Bigr),
\label{eq:BohnenblustSp}
\end{equation}
Here $\pi$ now runs through all unordered set partitions
$\mathcal{P}_n$ of $[n]$; by~$|\pi|$ we denote the number of
blocks in~$\pi$; and $m_i$ was the size of the particular
block~$\pi_i$. In the commutative case with weight $\theta=0$
we get the generalized integration by parts formula

\begin{equation}
    \sum_{\sigma\in S_n}
    R\Bigl( R\bigl( \cdots R(X_{\sigma_1})X_{\sigma_2} \cdots  \bigl) X_{\sigma_n} \Bigr)
    = \prod_{j=1}^{n} R\bigl(X_j\bigr).
\end{equation}
Also, for $n>0$ and $X_1 = \dots = X_n = X$ we find
\begin{equation}
    R\Bigl( R\bigl( \cdots R(X)X \cdots  \bigr) X \Bigr)
    = \frac{1}{n!}\sum_{\pi \in \mathcal{P}_n} (-\theta)^{n - |\pi|}
      \prod_{\pi_i \in \pi}(m_i-1)! \ R\bigl( X^{m_i}\bigr).
\end{equation}


\section{A new identity for Rota--Baxter algebras}

In
this section we provide a detailed proof of a Theorem announced
in~\cite{cefg2006,egp2007}.

First, recall that in a RB algebra $(A,R)$ we find for the RB
double product~(\ref{def:RBdouble}) $R(a \ast_\theta b)=R(a)R(b)$
which is just a reformulation of the Rota--Baxter relation. Next,
we introduce some notation.

Let $(A,R)$ be a RB algebra and $a_1,\ldots ,a_n$ be a collection
of elements in $A$. For any permutation $\sigma \in S_n$ we define
the element $T_\sigma(a_1,\ldots , a_n)$ as follows: define first
the subset $E_\sigma \subset\{1,\ldots ,n\}$ by $k \in E_\sigma$
if and only if $\sigma_{k+1}>\sigma_j$ for any $j\le k$. We write
$E_\sigma$ in the increasing order $1 \le k_1 < \cdots < k_p\le
n-1.$ Then we set:
\begin{equation}\label{tsigma}
    T_\sigma(a_1,\ldots ,a_n):=
    \Bigl(\cdots\bigl((a_{\sigma_1}\triangleright_\theta a_{\sigma_2})
    \triangleright_\theta \cdots\bigr)\triangleright_\theta a_{\sigma_{k_1}}\Bigr)
    \ast_\theta \cdots \ast_\theta
    \Bigl(\cdots\bigl((a_{\sigma_{k_p+1}}\triangleright_\theta a_{\sigma_{k_p+2}}) \triangleright_\theta
    \cdots\bigr)\triangleright_\theta a_{\sigma_{n}}\Bigr).
\end{equation}
There are $p+1$ packets separated by $p$ double RB products on the
right-hand side of the expression (\ref{tsigma}) above, and the
parentheses are set to the left inside each packet. Quite
symmetrically we define the element $U_\sigma(a_1,\ldots ,a_n)$ by
considering first the subset $F_\sigma\subset\{1,\ldots ,n\}$
defined by $l \in F_\sigma$ if and only if $\sigma_{l}<\sigma_j$
for any $j\ge l+1$. We write $F_\sigma$ in the increasing order:
$1\le l_1<\cdots < l_q\le n-1$. Then we set:
\begin{equation}\label{usigma}
    U_\sigma(a_1,\ldots ,a_n):=
    \Bigl(a_{\sigma_1}\triangleleft_\theta\bigl(\cdots (a_{\sigma_{l_1-1}}\triangleleft_\theta
    a_{\sigma_{l_1}})\bigr)\cdots\Bigr) \ast_\theta \cdots \ast_\theta \Bigl(a_{\sigma_{l_q+1}}
    \triangleleft_\theta\bigl(\cdots (a_{\sigma_{n-1}}\triangleleft_\theta a_{\sigma_n})\bigr)\cdots\Bigr).
\end{equation}
There are $q+1$ packets separated by $q$ double RB products on the
right-hand side of the expression (\ref{usigma}) above, and the
parentheses are set to the right inside each packet. The pre-Lie
operations $\triangleright_\theta$ and $\triangleleft_\theta$
involved in the right-hand side of equality (\ref{tsigma})
respectively (\ref{usigma}) are given by~(\ref{lem:leftRBpreLie})
respectively (\ref{lem:rightRBpreLie}).

\medskip

Following~\cite{lam1998} it is convenient to write a permutation
by putting a vertical bar after each element of $E_\sigma$ or
$F_\sigma$ according to the case. For example, for the permutation
$\sigma=(3261457)$ inside $S_7$ we have $E_\sigma=\{2,6\}$ and
$F_\sigma=\{4,5,6\}$. Putting the vertical bars:
\begin{equation}
\sigma=(32|6145|7),\hskip 20mm \sigma=(3261|4|5|7)
\end{equation}
we see that the corresponding elements in $A$ will then be:
 \allowdisplaybreaks{
\begin{eqnarray}
    T_\sigma(a_1,\ldots , a_7) &=&
    (a_3 \triangleright_\theta a_2)
    \ast_\theta \Bigl(\bigl((a_6\triangleright_\theta a_1)\triangleright_\theta a_4\bigr)
    \triangleright_\theta a_5\Bigr) \ast_\theta a_7,\\
    U_\sigma(a_1,\ldots, a_7)&=&
    \Bigl(a_3 \triangleleft_\theta\bigl(a_2\triangleleft_\theta (a_6\triangleleft_\theta a_1)\bigr)\Bigr)
    \ast_\theta a_4 \ast_\theta a_5 \ast_\theta a_7.
\end{eqnarray}}

\begin{thm} \label{newNCSpitzer} (New noncommutative Spitzer formula)
We have:
 \allowdisplaybreaks{
\begin{eqnarray}
    \sum_{\sigma\in S_n}
    R\Bigl( \cdots R\bigl(R(X_{\sigma_1})X_{\sigma_2}\bigr)\cdots X_{\sigma_{n}}\Bigr)
    &=&
    \sum_{\sigma\in S_n}
    R\bigl(T_\sigma(X_1,\ldots ,X_n)\bigr),
\label{eq:mainRB1}\\
    \sum_{\sigma\in S_n}
    R\Bigl(X_{\sigma_1} \cdots R\bigl(X_{\sigma_{n-1}}R(X_{\sigma_n})\bigr) \dots\Bigr)
    &=&\sum_{\sigma\in S_n}
    R\bigl(U_\sigma(X_1,\ldots,X_n)\bigr),
\label{eq:mainRB2}
\end{eqnarray}}
\end{thm}

In the weight $\theta=0$ case, the pre-Lie operations involved on
the right-hand side of the above identities reduce to $a
\triangleright_0 b=[R(a),b]=-b \triangleleft_0 a$. This case, in
the form (\ref{eq:mainRB2}), has been handled by C.S.~Lam
in~\cite{lam1998}, in the concrete situation when $A$ is a
function space on the real line, and when $R(f)$ is the primitive
of $f$ which vanishes at a fixed $T \in {\mathbb R}$. In the case
of a commutative RB algebra both identities agree since both the
left and right RB pre-Lie products~(\ref{lem:leftRBpreLie}),
(\ref{lem:rightRBpreLie}), respectively, agree. See \cite{KDF2007}
for analogous statements in the context of dendriform algebras.

\begin{proof}
The proof of (\ref{eq:mainRB1}) proceeds by induction on the
number $n$ of arguments, and (\ref{eq:mainRB2}) follows easily by
analogy. The case $n=2$ reduces to the identity:
\begin{equation}
    R\bigl(R(X_1) X_2\bigr) + R\bigl(R(X_2) X_1\bigr) = R(X_1)R(X_2) + R(X_2 \triangleright_\theta X_1),
\end{equation}
which immediately follows from the definitions. The case $n=3$ is
already non obvious and relies on considering the six permutations in
$S_3$:
\begin{equation*}
(1|2|3),\hskip 8mm (1|32),\hskip 8mm (2|31),\hskip 8mm
(21|3),\hskip 8mm (321),\hskip 8mm (312),
\end{equation*}
so that \allowdisplaybreaks{
\begin{eqnarray*}
    \sum_{\sigma\in S_3}
        R\Bigl(R\bigl(R(X_{\sigma_1})\, X_{\sigma_2}\bigr) X_{\sigma_3}\Bigr)
        &=&
         R(X_1)\, R(X_2)\, R(X_3)
       + R(X_1)\, R(X_3 \triangleright_\theta X_2)
       + R(X_2)\, R(X_3 \triangleright_\theta X_1) \nonumber \\
        & &
       + R(X_2 \triangleright_\theta X_1)\, R(X_3)
       + R\bigl((X_3 \triangleright_\theta X_2) \triangleright_\theta X_1\bigr)
       + R\bigl((X_3 \triangleright_\theta X_1) \triangleright_\theta X_2\bigr),
\end{eqnarray*}}
To prove the identity, we consider the following partition of the group $S_n$:
\begin{equation}\label{partition}
    S_n=S_n^n \amalg\coprod_{j,k=1}^{n-1}S_n^{j,k},
\end{equation}
where $S_n^n$ is the stabilizer of $n$ in $S_n$, and where
$S_n^{j,k}$ is the subset of those $\sigma\in S_n$ such that
$\sigma_j=n$ and $\sigma_{j+1}=k$. For $k \in \{1,\ldots,n-1\}$ we
set:
\begin{equation}
    S_n^k:=\coprod_{j=1}^{n-1}S_n^{j,k}.
\end{equation}
This is the subset of permutations in $S_n$ in which the two-terms
subsequence $(n,k)$ appears in some place. We have:
\begin{equation}
    S_n=\coprod_{k=1}^{n}S_n^{k}.
\end{equation}
Each $S_n^k$ is in bijective correspondence with $S_{n-1}$, in an
obvious way for $k=n$, and by considering the two-term subsequence
$(n,k)$ as a single letter for $k \not= n$. Precisely, in that
case, in the expansion of $\sigma\in S_n$ as a sequence $(\sigma_1,\ldots,\sigma_n)$, we replace the pair $(n,k)$ by $n-1$ and
any $j$, $k<j<n$ by $j-1$, so that, for example, $(2,1,5,3,4)\in
S_5^{3,3}$ is sent to $(2,1,4,3)$ by the bijection. For each
$\sigma \in S_n^{k}$ we denote by $\widetilde \sigma$ its
counterpart in $S_{n-1}$. Notice that for any $k \not=n$ and for
any $j \in \{1,\ldots ,n-1\}$, the correspondence $\sigma\mapsto
\sigmat$ sends $S_n^{j,k}$ onto the subset of $S_{n-1}$ formed by
the permutations $\tau$ such that $\tau_j=n-1$. The following
lemma is almost immediate:
\begin{lem}\label{bijections}
For $\sigma\in S_n^n$ we have:
\begin{equation}
    T_\sigma(a_1,\ldots ,a_n) = T_{\widetilde\sigma}(a_1,\ldots ,a_{n-1}) \ast_\theta a_n,
\end{equation}
and for $\sigma \in S_n^k, k<n$ we have:
\begin{equation}
    T_\sigma(a_1,\ldots ,a_n) = T_{\widetilde\sigma}(a_1,\ldots ,\widehat{a_k},\ldots, a_{n-1},
    a_n \triangleright_\theta a_k),
\end{equation}
where $a_k$ under the hat has been omitted.
\end{lem}

We rewrite the $n-1$-term sequence $(a_1,\ldots
,\widehat{a_k},\ldots, a_{n-1}, a_n \triangleright_\theta a_k)$ as
$(c_1^k,\ldots ,c_{n-1}^k)$. We are now ready to compute, using
the last lemma and the induction hypothesis:
 \allowdisplaybreaks{
\begin{eqnarray*}
 \lefteqn{
    \sum_{\sigma \in S_n} R\bigl(T_\sigma(a_1, \ldots ,a_n)\bigr)
     = \sum_{k=1}^n \sum_{\sigma \in S_n^k} R\bigl(T_\sigma(a_1,\ldots ,a_n)\bigr)}\\
    &=& \sum_{\tau \in S_{n-1}}
    R\Bigl(\bigl(R\bigl(\cdots R(R(a_{\tau_1}) a_{\tau_2}) \cdots \bigr) a_{\tau_{n-1}}\bigr) \ast_\theta a_n\Bigr)
    + \sum_{k=1}^{n-1} \sum_{\tau \in S_{n-1}}
    R\Bigl(R\bigl(\ldots R(R(c^k_{\tau_1}) c^k_{\tau_2}) \cdots \bigr) c^k_{\tau_{n-1}}\Bigr)\\
    &=& \sum_{\tau \in S_{n-1}}
    R\Bigl(R\bigl(R\bigl(\cdots R(R(a_{\tau_1}) a_{\tau_2}) \cdots \bigr) a_{\tau_{n-1}}\bigr) a_n\Bigr)
    - \sum_{\tau \in S_{n-1}}
    R\Bigl(R\bigl(\cdots R(R(a_{\tau_1}) a_{\tau_2}) \cdots \bigr) a_{\tau_{n-1}} \tilde{R}(a_n)\Bigr)\\
    & &+ \sum_{k=1}^{n-1} \sum_{\tau \in S_{n-1}}
    R\Bigl(R\bigl(\ldots R(R(c^k_{\tau_1}) c^k_{\tau_2}) \cdots (a_n\triangleright_\theta a_k) \cdots \bigr) c^k_{\tau_{n-1}}\Bigr)\\
    \end{eqnarray*}}
where $a_n\triangleright_\theta
a_k=R(a_n)a_k+a_k\tilde{R}(a_n)=c^k_{\tau_j}=c^k_{n-1}$ lies in
position $j$. Recall that $x \ast_\theta y=R(x)y-x\tilde{R}(y)$.
Using the definition of the pre-Lie operation
$\triangleright_\theta$ and the RB relation we get:
 \allowdisplaybreaks{
\begin{eqnarray*}
 \lefteqn{
    \sum_{\sigma\in S_n} R\bigl(T_\sigma(a_1,\ldots ,a_n)\bigr)}\\
    &=& \sum_{\tau \in S_{n-1}}
    R\Bigl(R\bigl(R\bigl(\cdots R(R(a_{\tau_1}) a_{\tau_2}) \cdots \bigr) a_{\tau_{n-1}}\bigr) a_n\Bigr)
    - \sum_{\tau \in S_{n-1}}
    R\Bigl(R\bigl(\cdots R(R(a_{\tau_1}) a_{\tau_2}) \cdots \bigr) a_{\tau_{n-1}} \tilde{R}(a_n)\Bigr)\\
    & &+ \sum_{k=1}^{n-1} \sum_{\tau \in S_{n-1} \atop \tau_1=n-1}
    R\Bigl(R\bigl(\cdots R\bigl(R(R(a_n) a_k) c_{\tau_2}^k\bigr) \cdots \bigr) c_{\tau_{n-1}}^k\Bigr)
    + \sum_{k=1}^{n-1} \sum_{\tau \in S_{n-1} \atop \tau_1=n-1}
    R\Bigl(R\bigl(\cdots R\bigl(R(a_k \tilde{R}(a_n)) c_{\tau_2}^k\bigr) \cdots \bigr) c_{\tau_{n-1}}^k\Bigr)\\
    & &+ \sum_{k=1}^n\sum_{j=2}^{n-1}\sum_{\tau \in S_{n-1} \atop \tau_j=n-1}
      R\biggl(R\Bigl(\cdots R\Bigl(R\Bigl(R\bigl( \cdots R(c_{\tau_1}^k) c_{\tau_2}^k\bigr) \cdots \Bigr) a_n\Bigr)
      a_k \cdots \Bigr) c_{\tau_{n-1}}^k\biggr)\\
    & &- \sum_{k=1}^n\sum_{j=2}^{n-1}\sum_{\tau\in S_{n-1} \atop \tau_j=n-1}
      R\biggl(R\Bigl(\cdots R\Bigl(R\bigl(\cdots R(c_{\tau_1}^k) c_{\tau_2}^k\bigr) \cdots  \tilde{R}(a_n)\Bigr)
       a_k\cdots \Bigr) c_{\tau_{n-1}}^k\biggr)\\
    & &+ \sum_{k=1}^{n-1}\sum_{j=2}^{n-1} \sum_{\tau\in S_{n-1} \atop \tau_j=n-1}
      R\biggl(R\Bigl(\cdots R\Bigl(R\Bigl(R\bigl(\cdots R(c_{\tau_1}^k) c_{\tau_2}^k\bigr) \cdots \Bigr)
      a_k \tilde{R}(a_n)\Bigr) \cdots \Bigr) c_{\tau_{n-1}}^k\biggr)
\end{eqnarray*}}
where $a_n$ lies in position $j$ (resp. $j+1$) in lines 4 and 5
(resp. in the last line) in the above computation, and where $a_k$
lies in position $j+1$ (resp. $j$) in lines 4 and 5 (resp. in the
last line). We can rewrite this going back to the permutation
group $S_n$ and using the partition (\ref{partition}):
 \allowdisplaybreaks{
\begin{eqnarray*}
    \lefteqn{\sum_{\sigma \in S_n} R\bigl(T_\sigma(a_1,\ldots ,a_n)\bigr)}\\
    &=& \sum_{\sigma\in S_n^n}
       R\Bigl(R\bigl(R\bigl(\cdots R(R(a_{\sigma_1}) a_{\sigma_2}) \cdots \bigr) a_{\sigma_{n-1}}\bigr) a_{\sigma_n}\Bigr)\\
    & & - \sum_{\sigma \in S_n^n}
       R\Bigl(R\bigl(\cdots R(R(a_{\sigma_1}) a_{\sigma_2}) \cdots\bigr) a_{\sigma_{n-1}}  \tilde{R}(a_{\sigma_{n}})\Bigr)\\
    & & +\sum\limits_{k=1}^{n-1} \sum\limits_{\sigma \in S_n^{1,k}}
       R\Bigl(R\bigl(R\bigl(\cdots R(R(a_{\sigma_1}) a_{\sigma_2}) \cdots \bigr) a_{\sigma_{n-1}}\bigr) a_{\sigma_n}\Bigr)\\
    & & +\sum\limits_{k=1}^{n-1} \sum\limits_{\sigma \in S_n^{1,k}}
       R\Bigl(\bigl(R\bigl(\cdots R( a_{\sigma_2}\tilde{R}(a_{\sigma_1})) \cdots \bigr)  a_{\sigma_{n-1}}\Bigr)\\
    & & +\sum_{k=1}^n\sum_{j=2}^{n-1}\sum_{\sigma \in S_n^{j,k}}
    R\biggl(R\Bigl( \cdots R\Bigl( R\Bigl( R\Bigl(\cdots R(a_{\sigma_1}) a_{\sigma_2} \cdots \Bigr) a_{\sigma_j}\Bigr)
    a_{\sigma_{j+1}}\Bigr)\cdots\Bigr) a_{\sigma_n}\biggr)\\
    & & - \sum_{k=1}^n\sum_{j=2}^{n-1}\sum_{\sigma \in S_n^{j,k}}
    R\biggl(R\Bigl( \cdots R\Bigl( R\Bigl(\cdots R(a_{\sigma_1}) a_{\sigma_2} \cdots
    \tilde{R}(a_{\sigma_j})\Bigr) a_{\sigma_{j+1}}\Bigr)\cdots\Bigr) a_{\sigma_n}\biggr)\\
    & & +\sum_{k=1}^{n-1}\sum_{j=2}^{n-1} \sum_{\sigma\in S_n^{j,k}}
    R\biggl(R\Bigl( \cdots R\Bigl( R\Bigl(\cdots R(a_{\sigma_1}) a_{\sigma_2} \cdots
    \bigl(a_{\sigma_{j+1}}\tilde{R}(a_{\sigma_j}) \bigr)\Bigr)\cdots\Bigr) a_{\sigma_n}\biggr).
\end{eqnarray*}}

Lines 2, 4 and 6 together give the left-hand side of
(\ref{eq:mainRB1}) whereas lines 3, 5, 7 and 8 cancel. More
precisely line 3 cancels with the partial sum corresponding to
$j=n-1$ in line 8, line 5 cancels with the partial sum
corresponding to $j=2$ in line 7, and (for $n\ge 4$), the partial
sum corresponding to some fixed $j\in\{3,\ldots ,n-1\}$ in line 7
cancels with the partial sum corresponding to $j-1$ in line 8.
This proves equality (\ref{eq:mainRB1}).
\end{proof}


\section{On the Magnus and Atkinson's recursions}
\label{ssect:ExpSolvAtkinson}

Now we return to Atkinson's recursions in Theorem~\ref{Atkinson1}.
We will focus only on the first equation in~(\ref{atkinsonEqs}).
Recall that the Spitzer algebra is naturally a graded Hopf
algebra. Moreover, we will assume the ---free--- Rota--Baxter
algebra to be of weight $\theta$. For computational convenience,
we consider the embedding of the Spitzer algebra $\mathcal{S}$
into $\mathcal{S}[[t]]$ defined on homogeneous elements $z$ of
degree $n$ in $\mathcal{S}$ by $z\longmapsto t^n\cdot z$. We agree
to identify $\mathcal{S}$ with its image in $\mathcal{S}[[t]]$, so
that this image is naturally provided with a graded Hopf algebra
structure and that the grading operation $Y$ now is naturally
given by $t\partial_t$. It is then obvious that the equation (a
generalized integral equation if we view the Rota--Baxter
operation $R$ as a generalized integral operator)
\begin{equation}
\label{eq:AtkinsonF}
    F = 1 + R(F \cdot Xt)
\end{equation}
is solved by the series $F=F(t):=\sum_{n \ge 0} t^n(RX)^{[n]}$
which is a group-like element in the Hopf algebra. The operation
of the Dynkin map on $F$ is given by
$$
    D(F(t))=F(t)^{-1} \cdot t\partial_t F(t) = \mathfrak{L}(t):= \sum\limits_{n > 0} t^n
    \mathfrak{L}_{\theta}^{(n)}(X),
$$
which, of course, implies the linear differential equation
$t\partial_t F(t)=F(t) \cdot \mathfrak{L}(t)$ and hence, by
comparing coefficients on both sides, the
recursion:
$$
   n(RX)^{[n]} = \sum_{k=0}^{n-1} (RX)^{[k]} \cdot \mathfrak{L}_{\theta}^{(n-k)}(X),
$$
which is a way to relate our Theorem~\ref{thm:ncBS1} to the
classical problem of finding explicit solutions to the first order
linear differential equations $\partial_tX=XA$. Now, the linear
differential equation $t\partial_tF(t)=F(t)\mathfrak{L}(t)$ is a
classical differential equation in noncommutative variables with
associated integral operator $P$ (the weight zero RB operator
$P=\int_0^t$) so that the equation can be solved with the usual
techniques for solving matricial or functional first order
differential equations. Actually, it is well-known that in the
noncommutative setting the differential equation $\partial_t
F(t)=F(t) \cdot \frac{1}{t}\mathfrak{L}(t)$ respectively the
integral equation:
$$
    F(t) = 1 + \int_0^t F(t') \cdot \frac{1}{t'}\mathfrak{L}(t')dt'
$$
can be solved via the exponential function. Notice that this
equation is a particular case of Atkinson's recursion with $F=1 +
P(F\cdot \hat{\mathfrak{L}})$, with
$\hat{\mathfrak{L}}(t)=\frac{\mathfrak{L}(t)}{t}$. Recall Magnus'
seminal work~\cite{magnus1954}. He proposed the exponential Ansatz
\begin{equation}
\label{eq:SolutionAtkinsonF}
    F(t) = \exp\bigl(\Omega[\hat{\mathfrak{L}}](t)\bigr),
\end{equation}
where $\Omega[\hat{\mathfrak{L}}](0)=0$.
Following~\cite{gelfand1995}, the series
for~$\Omega[\hat{\mathfrak{L}}]$
\begin{equation}
    \Omega[\hat{\mathfrak{L}}](t) = \sum_{n>0}\Omega_{n}t^n,
\label{eq:Magnus-Exp}
\end{equation}
can be expressed in terms of multiple integrals of nested
commutators. Magnus provided a differential equation which in turn
can be easily solved recursively for the terms $\Omega_{n}$
$$
    \frac{d}{dt}\Omega[\hat{\mathfrak{L}}](t)
        = \frac{-\hbox{ad}\, \Omega[\hat{\mathfrak{L}}]}
               {{\rm{e}}^{-\hbox{\eightrm ad}\,\Omega[\hat{\mathfrak{L}}]} - 1}(\hat{\mathfrak{L}})(t).
$$
which leads to the Magnus recursion
\begin{equation}
    \Omega[\hat{\mathfrak{L}}](t) = P\biggl(\hat{\mathfrak{L}} +
    \sum_{n>0} (-1)^n b_n \Big[\hbox{\rm ad} \bigl(\Omega[\hat{\mathfrak{L}}]\bigr)\Big]^{n} (\hat{\mathfrak{L}})\biggr)(t).
\label{Magnus}
\end{equation}
The coefficients are $b_n:=B_n/n!$ with~$B_n$ the Bernoulli
numbers. For $n=1,2,4$ we find $b_1=-1/2,b_2=1/12$ and
$b_4=-1/720$. We have $b_3=b_5=\cdots=0$.

Strichartz succeeded in giving a closed solution to Magnus'
expansion~\cite{strich1987}, see also \cite{miel1970,gelfand1995}.
He found
\begin{align}
    \Omega[\hat{\mathfrak{L}}](t) &= \sum_{n>0} \sum_{\sigma \in S_n}
    \frac{(-1)^{d(\sigma)}}{n^2\binom{n-1}{d(\sigma)\,}} \int_0^t\int_0^{t_1}\dots \int_0^{t_{n-1}}
    \big[[\dots[\hat{\mathfrak{L}}(t_{\sigma_n}),\hat{\mathfrak{L}}(t_{\sigma_{n-1}})]\dots],
    \hat{\mathfrak{L}}(t_{\sigma_1})\big]\,
    \,dt_n \dots dt_2\,dt_1.
\label{eq:StrichartzMagnus}
\end{align}
Here $d(\sigma)$ denotes the number of {\textit{descents}} in the
permutation $\sigma \in S_n$, that is, $d(\sigma )=|\{i<n,\ \sigma
(i)>\sigma (i+1)\}|$. In fact, more detail can be provided.
In~~\cite{gelfand1995} we find the following theorem for
$\Omega(t):=\Omega[\hat{\mathfrak{L}}](t)$.

\begin{thm} \cite{gelfand1995} The expansion of $\Omega(t)$ in terms
of the $\hat{\mathfrak{L}}_{\theta}^{(i)}$s writes
$$
 \Omega[\hat{\mathfrak{L}}](t) = \sum_{n>0} \sum_{\sigma \in S_n}
    \frac{(-1)^{d(\sigma)}}{n\binom{n-1}{d(\sigma)\,}} \int_0^t\int_0^{t_1}\dots \int_0^{t_{n-1}}
    \hat{\mathfrak{L}}(t_{\sigma_n})\cdots\hat{\mathfrak{L}}(t_{\sigma_2}) \hat{\mathfrak{L}}(t_{\sigma_1})
    \,dt_n \dots dt_2\,dt_1.
$$
\end{thm}

The coefficient of the term $\mathfrak{L}_{\theta}^{(i_1)} \cdots
\mathfrak{L}_\theta^{(i_m)}$ in the above expansion of $\Omega_k$
was also given in~\cite{gelfand1995}.
\begin{equation}
                   k \int_0^{t}dt_1 \int_0^{t_1}dt_2\dots \int_0^{t_{m-1}}dt_m
                           \sum_{\sigma \in S_m}
                           \frac{ (-1)^{d(\sigma)} }{ m \ \binom{m-1}{d(\sigma)} }
                           t^{i_{m}-1}_{\sigma_m} \dots
                           t^{i_{1}-1}_{\sigma_1}.
\end{equation}
In other terms, the Magnus expansion solves the Atkinson
recursion. We summarize our results in the following Theorem.

\begin{thm}\label{solution}
For an arbitrary weight $\theta$ RB algebra $A$, the Atkinson
recursion $F = 1+R(F \cdot Xt)$, $X \in A$ is solved by the
Strichartz' expansion:
$$
    F(t)=\exp \Bigg(\sum_{m>0}\sum\limits_{i_1+\cdots+i_n=m \atop i_1, \ldots, i_n>0}
             m \int_0^{t}dt_1 \int_0^{t_1}dt_2\dots \int_0^{t_{n-1}}dt_n
                           \sum_{\sigma \in S_n}
                           \frac{ (-1)^{d(\sigma)} }{ n \ \binom{n-1}{d(\sigma)} }
                           t^{i_{n}-1}_{\sigma_n} \dots
                           t^{i_{1}-1}_{\sigma_1}
                           \mathfrak{L}_{\theta}^{(i_1)}(X) \cdots
                           \mathfrak{L}_\theta^{(i_n)}(X)\Bigg)
$$
\end{thm}


\section{Solving Bogoliubov's counterterm recursion}
\label{ssect:Bogolibov}

Let us return to the Bogoliubov recursion~(\ref{BogoRec}), as
described in Section 3. As already noticed, its original setting is perturbative
quantum field theory (pQFT).
Connes and Kreimer associated to a renormalizable quantum field
theory the Hopf algebra $H=\bigoplus_{n=0}^{\infty}H_n$ of
ultraviolet (UV) superficially divergent one-particle irreducible
(1PI) Feynman graphs~\cite{ck2000,ck2001},
see~\cite{fg2005,ek2005,manchon2001} for reviews. This Hopf
algebra is polynomially generated by UV superficially divergent
1PI Feynman graphs, graded by the number of loops and
non-cocommutative. It is connected with the base field being, say,
the complex numbers $\mathbb{C}$. We exclude from considerations
theories with gauge symmetries, for which the Hopf algebra is
still commutative but, in general, is not polynomially
generated~\cite{kreimer2006,wvs2006a,wvs2006b} anymore.

The relevant quantities for the theory such as the Green functions
can be deduced in principle from Feynman rules --a prescription
associating to each (1PI) Feynman graph an integral. If the
integrals were convergent, Feynman rules would be imbedded into
the group $G(\mathbb{C}) \subset {\rm{Lin}}(H,\mathbb{C})$ of
algebra maps from $H$ to $\C$ with the Hopf algebra counit
$\varepsilon$ as the group unit. Here, ${\rm{Lin}}(H,\mathbb{C})$
denotes as usual the associative algebra of linear maps from $H$
to $\mathbb{C}$ equipped with the usual convolution product, $f
\star g:=m_{\mathbb{C}}(f \otimes g)\Delta$. Then, there would be
no obvious need for a renormalization process. However, these
integrals are most often divergent, hence can not be interpreted
as elements of $G(\mathbb{C})$, and require to be renormalized to
make sense.

The process of regularization, one of the most common ways to
proceed, is encoded in the change of the target space from
$\mathbb{C}$ to a commutative algebra $A$ which is supposed to be
equipped with an idempotent operator (called the renormalization scheme
operator) denoted by $R \in \End(A)$. Feynman rules
have, after regularization, a rigorous meaning as elements of
$G(A)$, see e.g. \cite{collins1984}. We denote the projections~$R(A):=A_-$ and
$\tilde{R}(A):=(1-R)(A)=A_+$ corresponding to the vector space splitting of
$A = R(A) \oplus \tilde{R}(A)$. In the dimensional regularization and minimal substraction
scheme, for example, $A=\mathbb{C}[\epsilon^{-1},\epsilon]]$,
$A_-=\epsilon^{-1}\mathbb{C}[\epsilon^{-1}]$,
$A_+=\mathbb{C}[[\epsilon]]$. As already pointed out,
one can show that for $(A,R)$ an idempotent Rota--Baxter algebra
such as $\mathbb{C}[\epsilon^{-1},\epsilon]]$, ${\rm{Lin}}(H,A)$
with the idempotent operator $\mathcal{R}$ defined by
$\mathcal{R}(f)=R \circ f$, for any $f \in {\rm{Lin}}(H,A)$, is a
noncommutative complete filtered unital Rota--Baxter algebra.

The Bogoliubov/Atkinson recursion allows then to decompose $G(A)$
as the set-theoretic product of its subgroups $G_-(A)$ and
$G_+(A)$:
$$
    \forall \gamma\in G(A),\ \exists !\gamma_-\in G_-(A),\ \gamma_+\in G_+(A),
     \ {\rm{ such\ that }}\ \gamma =\gamma_-^{-1}\ast \gamma_+,
$$
where we use the notations of Section 3.
Moreover $\gamma_+$ is a multiplicative map from the Hopf algebra
of Feynman diagrams $H$ to $\mathbb{C}[[\epsilon]]$, and
$\gamma^{ren}:=\lim_{\epsilon\rightarrow 0}\gamma_+$ is therefore
a well-defined element of $G(\mathbb{C})$, the ''renormalized
Feynman rule'' one was looking for to compute the relevant
properties of the theory.

Now, our results allow to give closed form expansions for
$\gamma_-$, $\gamma_+$ and $\gamma^{ren}$. Recall indeed from
Section 3 that $\gamma_-$ solves Atkinson's recursion:
 \allowdisplaybreaks{
\begin{eqnarray*}
    \gamma_- &=& e_A + \sum_{n>0} (\mathcal{R}a)^{[n]}
\end{eqnarray*}}
for $a:=e_A - \gamma$, $e_A:=\eta_A \circ \varepsilon$, and
analogously for $\gamma_+$ in terms of $\tilde{\mathcal{R}}$. In
conclusion, we get, in the weight $\theta:=-1$ RB algebra
$(Lin(H,A), \mathcal{R})$ the following theorems.

\begin{thm} \label{thm:BogoSol1}
We have for $a = e_A - \gamma $ \allowdisplaybreaks{
\begin{eqnarray}
\label{N}
    \gamma_- &=& e_A + \sum_{n>0} (\mathcal{R}a)^{[n]}\notag\\
             &=& e_A + \sum_{n>0} \sum\limits_{i_1+\cdots+i_k = n \atop  i_1,\dots,i_k > 0}
    \frac{\mathfrak{L}_{1}^{(i_1)}(a) \star \cdots \star \mathfrak{L}_{1}^{(i_k)}(a)}
         {i_1 (i_1 + i_2)\cdots(i_1+\cdots+i_k)}
\end{eqnarray}}
\end{thm}
where
$\mathfrak{l}_{1}^{(m)}(a)=(\mathfrak{l}_{1}^{(m-1)}(a)\triangleright_1
a)$ and $\mathfrak{L}_{1}^{(n+1)}(a):=R\bigl(\mathfrak{l}_{1}^{(n+1)}(a)\bigr)$ 
was defined in~(\ref{def:LamCtheta}). It is important to
underline that the expression given here applies in principle,
besides the minimal substraction and dimensional regularization
scheme, to any renormalization procedure which can be formulated
in terms of a Rota--Baxter structure.

The reader should notice the formal similarity of this solution
for $\gamma_-$ (the ``counterterm character'') with
Connes--Marcolli's formula for the universal singular frame,
see~\cite{cm2004,cm2004/2,cm2006}, see also
\cite{eunomia2006,manchon2001}. In fact, in the context of
dimensional regularization together with the minimal subtraction
scheme, there exists a linear map from the Connes--Kreimer Hopf
algebra of Feynman graphs to the complex numbers $\beta =
\sum_{n>0} \beta_n$, naturally associated to the counterterm
$\gamma_-$ (see \cite{eunomia2006} for details and a Lie theoretic
construction of $\beta$) and such that:
\begin{equation}
\label{eq:GammaBeta}
     \gamma_- = e_A+\sum_{n>0} \sum\limits_{k_1+\cdots+k_m = n \atop  k_1,\dots,k_m > 0}
     \frac{\beta_{k_1} \star \cdots \star \beta_{k_m}}
          {k_1(k_1 + k_2) \cdots (k_1 + \cdots + k_m)}\frac{1}{\epsilon^{n}}
\end{equation}
However, in spite of the similarity of formulas (\ref{N}) and
(\ref{eq:GammaBeta}) together with the techniques to obtain them,
the $\frac{\beta_i}{\epsilon^i}$'s do not coincide with the
$\mathfrak{L}_1^{(i)}(a)$'s. In fact, both are obtained from the action of the Dynkin operator $D=S\star Y$. They
are the homogeneous components of $D(\gamma_-)$ but, with respect to two
different graded Hopf algebra structures:
$\frac{\beta_i}{\epsilon^i}$ is simply the homogeneous component of $D(\gamma_-)$ in
the completed Hopf algebra ${\rm{Lin}}(H,A)$, whereas $\mathfrak{L}_1^{(i)}(a)$
is the image in ${\rm{Lin}}(H,A)$ of the homogeneous component of $D(\gamma_-)$ in
the completed Spitzer algebra ${\mathcal S}$ built on one generator (still written abusively $a$), with its associated Hopf structure as described in Section \ref{sect: RBandNCSF}.\\ 
From Theorem~\ref{solution} we conclude immediately:
\begin{thm} \label{thm:BogoSol2} We have
\allowdisplaybreaks{
\begin{eqnarray*}
    \gamma_-= \exp \biggl(
    \sum\limits_{m>0 \atop {i_1+\cdots+i_n=m \atop i_1,\ldots, i_n>0} }
    \!\!\!\!\int_0^{1}dt_1 \int_0^{t_1}dt_2 \! \dots \! \int_0^{t_{n-1}}dt_n
     \sum_{\sigma \in S_n}\
     \frac{ m(-1)^{d(\sigma)} }{ n \ \binom{n-1}{d(\sigma)} }\
     t^{i_{n}-1}_{\sigma_n} \dots t^{i_{1}-1}_{\sigma_1}\
     \mathfrak{L}_{1}^{(i_1)}(a) \star \cdots \star  \mathfrak{L}_{1}^{(i_n)}(a)\biggr)
\end{eqnarray*}}
\end{thm}

Returning to~\cite{gelfand1995} we may give a more combinatorial
expression for the formula in Theorem~\ref{thm:BogoSol2}, omitting
the dummy integrations. Recall the notion of a composition $I$ of
an integer $m$, i.e., a vector of positive integers, its parts,
$I:=(i_1,\ldots,i_k)$, of length $\ell(I):=k$ and weight
$|I|:=\sum_{j=1}^{k}i_j = m$. For instance, all compositions of
weight $3$ are $C_3:=\{(111),(21),(12),(3)\}$. The set of all
compositions $C:=\bigcup_{ m\ge 0}C_m$ is partially ordered by
reversed refinement, that is, $I \preceq J$ iff each part of $I$
is a sum of parts of $J$. We call $J$ finer than $I$. For
instance, $(1234) \preceq (11112211)$. Recall
$$
    \omega(I):=\bigl(i_1(i_1 + i_2)\cdots(i_1 + \cdots + i_k)\bigr)^{-1}.
$$
Let $J$ be a composition finer than $I$. Define
$\tilde{J}:=(J_1,\ldots,J_m)$ to be the unique decomposition of
the composition $J$, such that $|J_k|=i_k$ for $k=1,\ldots, m$.
Now define
$$
    \omega(J,I):= \prod_{k=1}^{m} \omega(J_k).
$$
We then deduce from~\cite{gelfand1995} (paragraphs 4.2 and 4.3)
the following formula for the exponent
$\Omega[\mathfrak{L}]:=\sum_{i>0} \Omega_i[\mathfrak{L}]$ in the
formula in Theorem~\ref{thm:BogoSol2}:
$$
    \Omega_n=\sum_{|J|=n}n\frac{(-1)^{\ell(J)-1}}{\ell(J)}\
    \sum_{J \preceq K =\{(k_1,\ldots ,k_p)\}}
    \omega(K,J)\
    \mathfrak{L}_{1}^{(k_1)}(a) \star \dots \star \mathfrak{L}_{1}^{(k_p)}(a).
$$

\newpage

\textbf{Acknowledgements}

\smallskip

The first named author acknowledges greatly the support by the
European Post-Doctoral Institute. He also thanks Laboratoire J. A.
Dieudonn\'e at Universit\'e de Nice Sophia-Antipolis and the
Institut for theoretical physics at Bielefeld University for warm
hospitality. The present work received support from the ANR grant
AHBE 05-42234.





\begin{thebibliography}{99}
\bibitem{bergeron2005}
    N. Bergeron and M. Zabrocki,
    ``The Hopf algebras of symmetric functions and quasisymmetric functions in
    non-commutative variables are free and cofree''
    \texttt{ArXiv:math.CO/059265}
\end{thebibliography}
\end{document}